\documentclass[12pt,a4paper]{amsart}

\usepackage[a4paper,left=2cm,right=2cm,top=15mm,bottom=15mm]{geometry}

\usepackage{amsmath}
\usepackage{amsthm}
\usepackage{amssymb}
\usepackage{amsfonts}
\usepackage{bm}
\usepackage{latexsym}
\usepackage{mathtools}
\usepackage{graphics}

\DeclareMathOperator{\intt}{int}

\DeclareMathOperator{\sign}{sign}

\usepackage[dvipsnames]{xcolor}

\definecolor{blue}{rgb}{0,0,0.8}
\definecolor{red}{rgb}{0.8,0,0}
\definecolor{darkgreen}{rgb}{0,0.6,0}

\newcommand{\dg}[1]{{\color{darkgreen}{#1}}}

\newcommand\restr[2]{{  \left.\kern-\nulldelimiterspace
  #1
  \vphantom{\big|}
  \right|_{#2}
  }}

\newcommand{\wx}{\widetilde{\mathbf{x}}}

\newcommand{\mol}{\overline{m}}
\newcommand{\mul}{\underline{m}}

\newcommand{\Pn}{{\mathcal P}_n}
\newcommand{\WW}{{\mathcal W}}

\newcommand{\I}{{\mathcal I}}

\newcommand{\W}{{\mathcal W}}

\newcommand{\JJ}{{\mathcal J}}

\newcommand{\Ji}{{\mathcal J}_i}
\newcommand{\Jj}{{\mathcal J}_j}

\newcommand{\Kij}{{\mathcal K}_{ij}}

\newcommand{\ai}{{\alpha}_i}

\newcommand{\RR}{{\mathbb R}}

\newcommand{\de}{{\delta}}

\newcommand{\al}{{\alpha}}

\newcommand{\ve}{{\varepsilon}}
\newcommand{\vk}{{\varepsilon}_k}

\newcommand{\vji}{{\varepsilon}_{j,i}}

\newcommand{\vki}{{\varepsilon}_{k,i}}

\newcommand{\yy}{\mathbf{y}}
\newcommand{\ee}{\mathbf{e}}
\newcommand{\zz}{\mathbf{z}}
\newcommand{\vv}{\mathbf{v}}
\newcommand{\uu}{\mathbf{u}}
\newcommand{\xx}{\mathbf{x}}
\newcommand{\ba}{\mathbf{a}}

\newcommand{\mm}{\mathbf{m}}
\newcommand{\ww}{\mathbf{w}}

\newcommand{\cc}{\mathbf{c}}

\newtheorem{theorem}{Theorem}
\newtheorem{corollary}{Corollary}
\newtheorem{lemma}{Lemma}
\newtheorem{proposition}{Proposition}
\newtheorem{definition}{Definition}

\theoremstyle{definition}

\newtheorem{remark}{Remark}

\begin{document}

\title[Weighted Lagrange interpolation]{Conjectures of Bernstein and Erd\H os \\ for weighted Lagrange interpolation \\ with Hermite weight}

\author{B\'ela Nagy}

\author{Szil\'ard Gy. R\'ev\'esz}\thanks{Sz. R. was supported in part by Hungarian National Foundation for Scientific Research, Grant \#'s K-146387, K-147153 and Excellence No. 151341.}

\author{Patr\'icia Szokol}
\thanks{P. Sz. was supported by the HUN-REN Hungarian Research Network through funding of her visiting researcher semester.}

\address{Béla Nagy
\newline \indent Bolyai Institute, University of Szeged,
\newline \indent Szeged, Aradi v\'ertanuk tere 1.
\newline \indent 6720 HUNGARY
}
\email{nbela@math.u-szeged.hu}

\address{Patrícia Szokol
\newline  \indent Faculty of Informatics, University of Debrecen  \newline \indent H-4002 Debrecen, \newline \indent P.O. Box 400,  HUNGARY \newline \indent and
\newline \indent  HUN-REN-UD Equations, Functions, Curves and their Applications Research Group} \email{szokol.patricia@inf.unideb.hu}

\address{Szil\'ard Gy. R\'ev\'esz
\newline  \indent HUN-REN A. R\'enyi Institute of Mathematics \newline \indent Budapest, Re\'altanoda utca 13-15. \newline \indent 1053 HUNGARY} \email{revesz.szilard@renyi.hu}

\date{\today}

\begin{abstract}
The classical Bernstein and Erd\H os Conjectures for algebraic Lagrange interpolation on a finite interval were settled in two consecutive papers by Kilgore and by de Boor and Pinkus in 1978. Extending these results to weighted interpolation on unbounded intervals has proved to be more difficult.

In a previous work, we established the Bernstein and Erd\H os Conjectures, together with the ``modified'' de Boor--Pinkus intertwining theorem, for weighted Lagrange interpolation on the half-line with exponential weights. Here we continue this line of research by considering the Hermite weight on the whole real line. Our main result is the Hermite-weighted analogue of the Bernstein and Erd\H os conjectures together with the appropriate version of the de Boor--Pinkus intertwining theorem.

A key ingredient in the classical proofs is the nonsingularity of certain derivative matrices associated with the local maxima of the Lebesgue function. This property already breaks down in the exponentially weighted setting, although the resulting difficulties can still be resolved. In the Hermite-weighted setting, however, the situation becomes more intricate, requiring a substantial refinement of the previous approach.
\end{abstract}

\maketitle



\section{Introduction}

The aim of the present paper is to establish an appropriate form of the Bernstein and Erd\H os Conjectures for Hermite-weighted polynomial interpolation on the real line. More precisely, we consider the Hermite-weight $w(t)=e^{-t^2}$, $t\in\RR$ and investigate the weighted Lagrange interpolation operator generated by a node system
$\xx\in S$
where $S:=\{(x_0,x_1,\ldots,x_n)\in\RR^{n+1}:\
x_0<x_1<\ldots<x_n\}$.

The original Bernstein Conjecture was formulated for algebraic polynomial interpolation on a compact interval more than a century ago \cite{Bernstein}. The problem remained open for several decades until it was eventually solved by Kilgore \cite{Kilgore}. Closely related questions were proposed by Erd\H os \cite{Erdos-1, Erdos-2}, and were subsequently settled through the work of Kilgore, Cheney, de Boor and Pinkus \cite{Kilgore-Cheney, CBoorPinkus}. In fact, de Boor and Pinkus proved a substantially stronger statement, often referred to as the intertwining property, which implies the Erd\H os Conjecture as a special case.

Following these achievements, considerable effort was devoted to extending the theory to weighted polynomial spaces and to more general interpolation settings; see, among others, \cite{Kilgore-84, Kilgore1985, Kilgore-87, K-AMH}. A common feature of essentially all previous approaches is their reliance on the nonsingularity of certain derivative matrices built from the local maxima functions of the corresponding Lebesgue function. The decisive role of this property was later clarified by Shi \cite{Shi}, who showed that, together with a mild properness assumption, nonsingularity already implies not only the Bernstein and Erd\H os Conjectures, but even the stronger intertwining property.

For exponentially weighted interpolation on the halfline, however, this classical scheme breaks down. It was shown in \cite{Exp} that interval maxima may occur at endpoints of the corresponding intervals, causing entire columns of the derivative matrices to vanish. Consequently, the standard nonsingularity argument cannot be repaired in that setting. Nevertheless, the Bernstein and Erd\H os Conjectures were established in \cite{Exp} by a completely different approach based on global topological considerations and a careful analysis of the structure of the local maxima.

The present paper continues this line of investigation in the Hermite-weighted setting. Although the Hermite weight shares several qualitative features with the exponential weight, the geometry of the interpolation problem is substantially different. The exponential case is formulated on a halfline and possesses only one outer interval. In contrast, the Hermite setting is naturally posed on the whole real line and therefore contains two outer intervals.

This distinction has important consequences for the derivative matrices. In the exponential case considered in \cite{Exp}, singularity arises when the maximum of the last interval polynomial is attained at the endpoint of the outer interval, leading to a vanishing column in the derivative matrix. In the Hermite case the phenomenon may occur simultaneously at both ends of the real line. As a consequence, not merely one but potentially two columns of the derivative matrix may vanish identically. Equivalently, after transposition, one may encounter two zero rows. Thus the failure of the nonsingularity framework is even more pronounced than in the exponential case.

Therefore, the Hermite-weighted problem provides a natural and significantly more demanding test case for methods that do not rely on nonsingularity. The purpose of the present work is to show that despite these difficulties, the Bernstein and Erd\H os Conjectures remain valid in this setting.

As in the exponential case, we consider the weighted space $$C_w(\RR):=\left\{f\in C(\RR): \|f\|_w:=\|fw\|_{\infty}<\infty \right\},$$ together with the weighted polynomial space $\W_n:=w\Pn$.
 In the present paper,
however, the weight is mostly the Hermite weight $w(t)=e^{-t^2}$; for some more general statements
we will also use positive, continuous weights satisfying $\lim_{\pm \infty} w(t)t^n =0$, subject to appropriate further conditions. 

For a node system
$\xx=(x_0,\ldots,x_n)\in S$ 
there are unique weighted fundamental functions
\begin{equation}\label{def:h_k}
h_k(\xx;t)=\frac{w(t)}{w(x_k)}\prod_{\substack{j=0\\ j\neq k}}^n\frac{t-x_j}{x_k-x_j}, \qquad k=0,\ldots,n,
\end{equation}
satisfying $h_k(x_j)=\delta_{jk}$
where $\delta_{jk}$ is the Kronecker delta. 

The corresponding weighted Lebesgue function is given by
\begin{equation}
L_n(w,\xx,t)=\sum_{k=0}^n |h_k(\xx,t)|=\sum_{k=0}^{n} \vk(t) h_k(t)
\end{equation}
where
$$
\vk(t):=\sign(h_k(\xx,t)) = \begin{cases} (-1)^{k-i} \qquad & \textrm{if} \quad t\in \intt I_i, ~ 1 \le  i \le k \le n
\\ (-1)^{k+1-i} & \textrm{if} \quad t\in \intt I_i, ~0 \le  k < i \le n+1
\\ \delta_{ik} & \textrm{if} \quad t=x_i, \quad i=0,1,\ldots,n
\end{cases},
$$
and, for a fixed node system
$\xx=(x_0,\ldots,x_n)\in S$ 
the nodes partition the real line into the intervals
\[
I_0:=(-\infty,x_0],\qquad
I_i:=[x_{i-1},x_i]\quad (1\le i\le n),
\qquad
I_{n+1}:=[x_n,\infty).
\]

On each interval the signs of the weighted fundamental functions remain constant. Therefore, restricting the Lebesgue function to a given interval $I_i$, we obtain a weighted polynomial

\begin{equation}\label{def:Pi}
P_i(t):=P_i(\xx,t):=L_n(w,\xx,t)\big|_{I_i}=\sum_{k=0}^{n}\ve_{k,i}h_k(\xx,t),\qquad i=0,\ldots,n+1,
\end{equation}
where
\begin{equation}\label{epsilonki}
\vki=(-1)^{k+1-i+\chi_{i\le k}}\quad (0\le k \le n, ~0 \le i \le n+1).
\end{equation}

Then, we have that $P_i(x_{i-1})=P_i(x_{i})=1$, ($i=1,\ldots,n$), also $P_0(x_0)=1$ and $P_{n+1}(x_n)=1$. Moreover, $P_0=(-1)^nP_{n+1}$. The norm of the interpolation operator is therefore determined by the local maxima $$m_i(\xx):=\max_{t\in I_i}P_i(\xx,t), \qquad i=0,\ldots,n+1,$$ since

\[
\|\mathcal L_n\|_w=\max_{t\in\RR}L_n(w,\xx,t)=\max_{0\le i\le n+1}m_i(\xx).
\]

Consequently, the problem of finding optimal interpolation nodes, that is, node systems minimizing the norm of the interpolation operator, is equivalent to minimizing the maximum of the quantities $m_i(\xx)$.

For a more detailed discussion of the weighted Lebesgue function and the related notation, we refer the reader to \cite{Exp}. The setting considered there differs from the present one in that interpolation is studied on the half-line $[0,\infty)$ with an exponential weight, whereas here we work on $\RR$ with the Hermite weight.

The Bernstein Conjecture predicts that an optimal node system $\xx$ is characterized by the equioscillation property

\[
m_0(\xx)=m_1(\xx)=\cdots=m_{n+1}(\xx).
\]

The corresponding Erd\H os Conjecture asserts that if $\xx$ is an optimal node system and $\yy\neq \xx$  is any other node system, then the interval maxima of \(\yy\) must straddle the common optimal value. More precisely, there exist indices $i$ and $j$ such that

\[
m_i(\yy)<m_i(\xx)
\qquad\text{and}\qquad
m_j(\yy)>m_j(\xx).
\]

In other words, no competing node system can have all interval maxima simultaneously above or simultaneously below the equioscillating optimal level.
The stronger intertwining property of de Boor and Pinkus \cite{CBoorPinkus} states that the same conclusion holds for any two distinct node systems $\xx\neq\yy$, without assuming that either of them is optimal.

Our main result shows that the Bernstein and Erd\H os assertions remain valid for Hermite-weighted interpolation. A crucial ingredient in the classical proofs is the nonsingularity of the $(n+1)\times(n+1)$ submatrices $A_k$, obtained by deleting the $k$th column from
the transpose of
the Jacobian matrix, i.e., from the matrix
\begin{equation}\label{def_matrix_A}
A=\left(
\frac{\partial m_i}{\partial x_j}(\xx)
\right)_{\substack{j=0,\ i=0}}^{\substack{n,\ n+1}}.
\end{equation}

In \cite{K-AMH}, Kilgore claimed that these matrices remain nonsingular not only in the classical algebraic case, but also for exponentially and Hermite-weighted interpolation, and consequently formulated the corresponding Bernstein and Erd\H os theorems in these settings. However, the claimed nonsingularity of the matrices $A_k$ was disproved in \cite{SzAMH} by explicit counterexamples (see Examples~1 and ~3 of \cite{SzAMH} for the exponential and Hermite cases, respectively). As a consequence, Kilgore's exponentially- and Hermite weighted versions of the Bernstein and Erd\H os theorems reverted to conjectures.

In our previous work \cite{Exp}, we proved the exponentially weighted case without assuming the nonsingularity of the matrices $A_k$. Nevertheless, even in that setting substantial new ideas were required to compensate for the failure of the classical argument.

The Hermite-weighted setting turns out to be somewhat more delicate. Unlike the exponential half-line problem, the real line possesses two outer intervals, $I_0=(-\infty,x_0], I_{n+1}=[x_n,\infty)$,
which are not bounded by interpolation nodes on both sides. Consequently, the maxima of the
outer 
weighted polynomials $P_0$ and $P_{n+1}$, defining $m_0(\xx)$ and $m_{n+1}(\xx)$, need not be attained in the interiors of $I_0$ and $I_{n+1}$, but may instead occur at the boundary nodes $x_0$ and $x_n$, respectively.

This phenomenon leads to a substantially stronger form of degeneracy than in the exponential case. Indeed, one or even both of the columns of $A$ corresponding to the ``outer maxima'' (i.e., $\nabla m_0(\xx)$ and $\nabla m_{n+1}(\xx)$) may vanish in the derivative matrices appearing in the Bernstein--Erd\H os approach, whereas in the exponential setting only one such degeneracy can occur.

The approach developed in the present paper is also related to several earlier methods designed to replace the classical nonsingularity argument. In particular, de Boor and Pinkus regarded the homeomorphism property of an appropriate difference mapping as one of the central conclusions of their work \cite{CBoorPinkus}. A related formulation was later investigated by Shi \cite{Shi}. In the absence of nonsingularity, this homeomorphism property provides a natural substitute for the classical Jacobian argument, and it will play a central role in our proof as well. Our use of this idea is also motivated by its successful application in other approximation-theoretic settings; see, for example, \cite{TLMS2018,Homeo}. Nevertheless, before applying the homeomorphism approach, we need to make substantial detours, because here even two columns of the matrix $A$ may be vanish, whence the rank of the matrix can be lower than $n$, excluding a  homeomorphic mapping of node systems in general.

A second source of motivation comes from Fenton's equioscillation theorem for sums of translates of concave kernels \cite{Fenton}. Building on this idea, and with the classical results of Kilgore and de Boor--Pinkus in mind, general minimax and equioscillation theorems were obtained in several settings where nonsingularity was unavailable \cite{Ural,Sbornik,JMAA}. Related periodic results go back to Kendall, Hardin and Saff \cite{HardinKendallSaff} and were further developed in \cite{TLMS2018}. In those problems, however, the failure of nonsingularity could be compensated for by additional structural properties, most notably the concavity of the underlying kernels and the fact that the summands depended only on a single difference variable. These features are absent in the present Hermite-weighted interpolation problem, and hence the corresponding one-dimensional techniques cannot be applied directly.

\section{The regular set and its boundary structure}

\noindent
We now describe the global oscillatory structure of the ``outer polynomials\footnote{The terminology refers to the fact that these polynomials belong to intervals lying outside the node points.}'' $P_{n+1}(\xx,\cdot)$ and $P_0(\xx,\cdot)$. Following the terminology of Milev and Naidenov \cite{Milev-Naidenov-JAT-2006}, an oscillating polynomial (that is, a polynomial whose zeros are all real, simple, and as numerous as its degree), is said to satisfy Property (P) with parameters $\delta_0,\delta_n\in\{0,1\}$, if its derivative has exactly one simple zero between every two consecutive zeros of the polynomial, together with $\delta_0$ additional simple zero to the left of the smallest zero and $\delta_n$ additional simple zero to the right of the largest zero.

For the Hermite-weighted setting on $\RR$, Milev and Naidenov (see \cite{Milev-Naidenov-JAT-2006}) proved that Property (P) holds for oscillating polynomials with $\delta_0=\delta_n=1$.

As will be shown later in  Section \ref{sec:zeros},
for every $\xx\in S$ the weighted polynomial $P_{n+1}(\xx,\cdot)$ is oscillating and therefore has exactly $n$ simple zeros. Moreover, since $P_0=(-1)^nP_{n+1}$, the zeros of the two weighted polynomials coincide. Thus, if
\[
y_1<\ldots<y_n
\]
denote the zeros of $P_{n+1}(\xx,\cdot)$, then these are also the zeros of $P_0(\xx,\cdot)$.

Applying the Rolle theorem the derivative functions $(P_0)'_t$ and $(P_{n+1})'_t$ have exactly $n-1$ simple interior zeros
\[
 y_1 < w_1 < y_2 < \ldots < w_{n-1} < y_n,
\]
and the property (P) guarantees, that there exist two further zeros $w_0$ and $w_n$ of $(P_0)'_t(\xx,\cdot)$ and $(P_{n+1})'_t(\xx,\cdot)$ such that $w_0\in (-\infty, y_1)$ and $w_{n}\in (y_{n},\infty)$.

Since $P_{n+1}(\xx,x_n)=1$, $\lim_{t\to \infty}P_{n+1}(\xx,t)=0$ and $P_{n+1}(\xx,\cdot)$ has only one critical point to the right of its largest zero, namely $w_n$, it follows that $P_{n+1}(\xx,\cdot)$ is strictly increasing on $(y_n,w_n)$ and strictly decreasing on $(w_n,\infty)$.

Consequently, the maximum of $P_{n+1}(\xx,\cdot)$ on the outer interval $[x_n,\infty)$ is attained at a unique point, which we denote by $z_{n+1}:=z_{n+1}(\xx)$; in fact,
$z_{n+1}=\max(w_n,x_n)$.
Similarly, the maximum of $P_0(\xx,\cdot)$ on $(-\infty,x_0]$ is attained at a unique point, denoted by
$z_0:=z_0(\xx)$. More generally, throughout the paper, $z_i:=z_i(\xx)$ will denote the unique\footnote{The uniqueness of $z_i(\xx)$ will be established in Section~4.} point at which $P_i(\xx,\cdot)$ attains its maximum on the interval $I_i$.

The present section is the Hermite-weighted analogue of the corresponding section in \cite{Exp}. As in the exponential setting, the arguments are carried out in the framework of Extended Chebyshev--Haar Systems (ECHS)\footnote{Both names occur in the literature. Although the concept originated with Chebyshev, it first appeared in print in Markov's work \cite{Markov} and was later introduced independently by Haar \cite{Haar}.}. Recall that a vector space of functions is called a rank $n$ ECHS, if its dimension is $n$ and the unique function from the family having $n$ zeros -- counting according to multiplicities -- is the identically zero function. Obviously, the weighted polynomials $P_0,P_1,\ldots,P_{n+1}$ belong to an ECHS of rank $n+1$.

The following identity originates in the work of Kilgore and Cheney \cite{Kilgore-Cheney}, who established a related formula for the neighboring weighted polynomials $P_{j-1}$ and $P_j$ under the additional assumption ${\bf 1}\in\W_n$. In Lemma~2 of \cite{Exp} this result was extended to general ECHS, without requiring the constant function to belong to the system. We recall it here, as it will be used repeatedly later, especially in Section~\ref{sec:linearindependence}.

\begin{lemma}\label{l:magic_0}
Let $\W_n\subset C^1(\RR)$ be a rank $n+1$ ECHS of continuously differentiable functions, and let $\xx\in S$. Then, for every $j=0,\ldots,n$ and $i=0,\ldots,n+1$,
\begin{equation}\label{magic-middleform_0}
\frac{\partial P_i}{\partial x_j}(\xx,t)
=
-h_j(\xx,t)\,(P_i)'_t(\xx,x_j),
\qquad t\in\RR.
\end{equation}
\end{lemma}

The present section is the Hermite-weighted analogue of Section~2 in \cite{Exp}. The underlying ECHS framework remains unchanged, and the analysis of the right outer interval $I_{n+1}=[x_n,\infty)$ is identical to that in \cite{Exp}. The essential novelty is the presence of the other outer interval $I_0=(-\infty,x_0]$, where the function $m_0$ exhibits the same type of boundary behaviour as $m_{n+1}$. Consequently, every argument concerning the right outer interval $I_{n+1}$ has a symmetric counterpart for the left outer interval $I_0$. Throughout this section we therefore refer to the corresponding results of \cite{Exp} whenever the analysis is unchanged, and present only the additional arguments required for the left endpoint.

We introduce the \emph{regular set}, that can be considered as a nice part of $S$.

\begin{definition}\label{def:X_and_complement}
Let $\WW\subset C^1(\RR)$ be a rank $n+1$ ECHS. The \emph{regular set} $X\subset S$ is defined as
\[
X:=\{\xx\in S:\ m_0(\xx)>1,\ m_{n+1}(\xx)>1\},
\]
and concerning the (relative) complement $X^c:=S\setminus X$ of the regular set $X$, we introduce
\begin{align*}
X_-:=\{\xx\in S:\ m_0(\xx)=1\},
\qquad
X_+:=\{\xx\in S:\ m_{n+1}(\xx)=1\}.
\end{align*}
\end{definition}
Then the complement of the regular set can be written as the non-disjoint union
\[
X^c = X_- \cup X_+.
\]
We note that $S\setminus X_-=\{\xx\in S:\ m_0(\xx)>1\}$, respectively $S\setminus X_+=\{\xx\in S:\ m_{n+1}(\xx)> 1\}$, and therefore both sets, as well as $X$ are open by the continuity of $m_0$ and $m_{n+1}$. Correspondingly, $X_+, X_-, X^c$ are relatively closed in $S$.

In contrast with the classical situation treated by Kilgore and Cheney \cite{Kilgore-Cheney}, the outer interval maxima need not be attained in the interior of the corresponding intervals. Node systems outside the regular set are precisely those configurations for which at least one of the two exterior maxima is attained at the corresponding endpoint of the intervals $I_0$ and $I_{n+1}$ (i.e., at $x_0$ or  at $x_n$).

As in the exponential case (see Proposition 1 in \cite{Exp}), we can also characterize the elements of X (more precisely, those of $S\setminus X_+$ and $S\setminus X_-$).

\begin{proposition}\label{prop:equiv-hermite}
Let $\WW\subset C^1(\RR)$ be a rank $n+1$ ECHS satisfying Property (P) with $\de_0=\de_n=1$. Then, for every $\xx\in S$, the following are equivalent:

\begin{enumerate}
\item[1)] $m_{n+1}(\xx)>1$, i.e., $\xx\in S\setminus X_+$;

\item[2)] the unique zero $w_{n}(\xx)$ of $(P_{n+1})'_t(\xx,\cdot)$ on $(y_n,\infty)$ satisfies $w_{n}(\xx)>x_n$;

\item[3)] the unique
maximizer 
of $P_{n+1}(\xx,\cdot)$ on $[x_n,\infty)$ lies in $(x_n,\infty)$;

\item[4)] $z_{n+1}=w_{n}>x_n$;

\item[5)] $(P_{n+1})'_t(\xx,x_n)>0$.
\end{enumerate}

Similarly, the following are equivalent:

\begin{enumerate}
\item[6)] $m_0(\xx)>1$, i.e., $\xx\in S\setminus X_-$;

\item[7)] the unique zero $w_0(\xx)$ of $(P_{0})'_t(\xx,\cdot)$ on $(-\infty,x_0]$ lies in $(-\infty,x_0)$;

\item[8)] the unique maximiser $z_0$ of $P_0(\xx,\cdot)$ on $(-\infty,x_0]$ lies in $(-\infty,x_0)$;

\item[9)] $z_{0}=w_{0}<x_0$;

\item[10)] $(P_0)'_t(\xx,x_0)<0$.
\end{enumerate}
\end{proposition}

\begin{proof}
The characterization of $S\setminus X_+$ follows exactly the same argument as Proposition~1 of \cite{Exp}, applied to the right outer interval $[x_n,\infty)$. In particular, the proof relies on Lemma~\ref{l:magic_0}, which is the analogue of the corresponding identity used in \cite{Exp}. The characterization of $S\setminus X_-$ is then obtained by the symmetric argument on the left outer interval $(-\infty,x_0]$.
\end{proof}

We denote the parameter spaces
\[
S_0:=\{(x_0,\ldots,x_{n-1})\in \RR^n: x_0<\cdots<x_{n-1}\}, \quad
S_n:=\{(x_1,\ldots,x_n)\in \RR^n:x_1<\cdots<x_n\}.
\]

The above characterization immediately leads to the corresponding threshold values associated with the two outer intervals.

\begin{lemma}\label{l:xstar-hermite}
Let $\WW\subset C^1(\RR)$ be a rank $n+1$ ECHS satisfying Property (P) with $\de_0=\de_n=1$.
Fix $\wx=(x_1,\dots,x_{n-1})$ and $x_0$ such that $(x_0,\wx)\in S_0$, i.e., $x_0<x_1$.

Then there exists a unique number $x_n^*>x_{n-1}$ such that
\begin{equation}
(P_{n+1})'_t((x_0,\wx,x_n^*),x_n^*)=0.
\label{d:xnstar}
\end{equation}

Moreover, if $x_n>x_n^*$, then the node system $\xx=(x_0,\wx,x_n)$ satisfies that $(P_{n+1})'_t(\xx,x_n)<0$, hence
$\xx \in\intt (X_+)$,
and on the other hand, if $ x_n<x_n^*$,  then $(P_{n+1})'_t(\xx,x_n)>0$, so $\xx\in S\setminus X_+$.

Similarly, fixing $(\wx,x_n)\in S_n$, there exists a unique number $x_0^*<x_1$ such that
\begin{equation}
(P_0)'_t((x_0^*,\wx,x_n),x_0^*)=0.
\label{d:xnzero}
\end{equation}

Moreover, if $x_0<x_0^*$, then the node system $\xx=(x_0,\wx,x_n)$ satisfies that $(P_0)'_t(\xx,x_0)>0$, hence
$\xx \in\intt(X_-)$,
and on the other hand, if $x_0>x_0^*$, then $(P_0)'_t(\xx,x_0)<0$,
so $\xx\in S\setminus X_-$.
\end{lemma}

\begin{proof}
First, fix $(x_0,\wx)\in S_0$. The argument proving Lemma~3 of \cite{Exp} applies without modification
to the right outer interval $I_{n+1}=[x_n,\infty)$, yielding a unique threshold value $x_n^*>x_{n-1}$ satisfying $(P_{n+1})'_t((x_0,\wx,x_n^*),x_n^*)=0$.
Moreover,
\[
x_n>x_n^* \quad\Longrightarrow\quad (P_{n+1})'_t(\xx,x_n)<0, \quad \text{ while } \quad
x_n<x_n^* \quad\Longrightarrow\quad (P_{n+1})'_t(\xx,x_n)>0.
\]
Consequently, by Proposition~\ref{prop:equiv-hermite},
\[
x_n>x_n^*\Longrightarrow \xx\in \operatorname{int}(X^c), \quad \text{ whereas } \quad
x_n<x_n^*\Longrightarrow\xx\in S\setminus X_+.
\]
Thus every half-line
\[
\ell^+_{(x_0,\wx)}=\{(x_0,\wx,t):t>x_{n-1}\}
\]
contains exactly one threshold point $(x_0,\wx,x_n^*)$, and this point separates the region $S\setminus X_+$ from $\intt (X^c)$ along that half-line.

The same reasoning, applied to the interval $I_0$, yields the existence and uniqueness of $x_0^*$, together with the corresponding sign characterization.
\end{proof}

\begin{proposition}\label{prop:fat-hermite}
Assume that $\WW\subset C^1(\RR)$ is a rank $n+1$ ECHS satisfying Property (P) with $\de_0=\de_n=1$.
Then, there exist continuous functions $\phi_+: S_0\to\RR$, and $\phi_-: S_n\to\RR$,
such that $\phi_+(x_0,\wx)=x_n^*$ and $\phi_-(\wx,x_n)=x_0^*$, where $x_n^*$ and $x_0^*$ are the unique threshold values from Lemma~\ref{l:xstar-hermite}
defined in \eqref{d:xnstar} and
\eqref{d:xnzero}.

Furthermore, the sets $X_-$ and $X_+$ are \emph{fat}, i.e.,
\[
X_\pm=\overline{\intt(X_\pm)}.
\]
Consequently, $X^c=X_-\cup X_+$ is fat as well.
\end{proposition}

\begin{proof}
Define according to \eqref{d:xnstar} and
\eqref{d:xnzero}
\[
\phi_+(x_0,\wx):=x_n^*,
\qquad
\phi_-(\wx,x_n):=x_0^*.
\]
The continuity of both functions follows exactly as in Proposition~2 of \cite{Exp}. Indeed, the only ingredients required there are the continuity of the corresponding derivative function, the uniqueness of the threshold value, and the sign change established in Lemma~\ref{l:xstar-hermite}.

We prove the fatness only for $X_-$, the argument for $X_+$ being completely symmetric. By Lemma~\ref{l:xstar-hermite},
\[
X_-=\{\xx\in S:\ x_0\le \phi_-(\wx,x_n)\}.
\]
Hence every relative boundary point of $X_-$ satisfies
\[
x_0=\phi_-(\wx,x_n).
\]
Again by Lemma~\ref{l:xstar-hermite}, every such point can be approximated by points from $\intt(X_-)$ obtained by decreasing  the first node $x_0$. Therefore, $X_-=\overline{\intt(X_-)}$.

The proof for $X_+$ is completely analogous. Finally, $X^c=X_-\cup X_+$ is fat as the union of two fat sets.
\end{proof}

\section{Location of zeros of $P_i$}\label{sec:zeros}

In this section we analyse the location of the zeros of the interval polynomials $P_i$.
We first determine the number and location of the zeros of the functions $P_i$, $i=0,\ldots,n+1$. The analysis is based only on the sign pattern of the functions $P_i$ at the interpolation nodes and on the fact that the corresponding weighted polynomial system is an ECHS.

As in the case of exponentially weighted interpolation, (see \cite{Exp}), we occasionally pass to the polynomial parts
\[
p_i(\xx,t):=\frac{P_i(\xx,t)}{w(t)}, \qquad i=0,\ldots,n+1.
\]
Since $w(t)>0$ on $\RR$, the zeros of $p_i$ and $P_i$ coincide. Hence, when we speak about the degree of $P_i$, we mean the degree of $p_i$.

For $i=1,\ldots,n$, the function $P_i$ takes alternating signs at the nodes $x_0,\ldots,x_n$, except at the two neighbouring nodes $x_{i-1}$ and $x_i$, where $P_i(\xx,x_{i-1})=P_i(\xx,x_i)=1$.
Consequently, $P_i$ has one zero in each interval $I_k=[x_{k-1},x_k]$ $k=1,\ldots,n,\quad k\neq i$, and no zero in $I_i$. These zeros are simple and unique in their respective intervals. Indeed, if one of them were multiple zero, or if there were an additional zero in any of these intervals, then $P_i$ would have at least $n+1$ zeros counted with multiplicities, which is
impossible in an ECHS of dimension rank $n+1$.

The two outer interval polynomials $P_0=(-1)^nP_{n+1}$ are different. They change sign at every pair of consecutive nodes. Therefore $P_0=(-1)^nP_{n+1}$ has exactly one simple zero in every interval $I_k=[x_{k-1},x_k],$ $k=1,\ldots,n$.

It remains to decide whether the interior polynomials $P_i$, $i=1,\ldots,n$, may have an additional zero in one of the outer intervals $I_0=(-\infty,x_0]$ and $I_{n+1}=[x_n,\infty)$.

This is determined by the leading coefficient of the polynomial part $p_i$.

Similarly to the exponential case, let $a_i$ denote the coefficient of $t^n$ in $p_i$, and let $b_k$ be the leading coefficient of the classical Lagrange fundamental polynomial $\ell_k$. Then
\[
b_k=\prod_{\substack{j=0\\j\ne k}}^n\frac{1}{x_k-x_j}=(-1)^{n-k}|b_k|,
\qquad k=0,\ldots,n.
\]
Consequently,
\begin{equation}\label{ai}
\begin{aligned}
a_i= \sum_{k=0}^n \frac{\vki}{w(x_k)} b_k &= \sum_{k=0}^n \frac{\vki}{w(x_k)} (-1)^{n-k} |b_k|
= (-1)^{n+1+i} \left\{\sum_{k=0}^{i-1} \frac{|b_k|}{w(x_k)} - \sum_{k=i}^{n} \frac{|b_k|}{w(x_k)} \right\}.
\end{aligned}
\end{equation}

From this it is clear that $a_{n+1}>\cdots > (-1)^{n+1+i}a_i > \cdots > (-1)^{n+1} a_0$. Moreover, by \eqref{ai} we get that $a_{n+1}>0$ and $P_0=(-1)^nP_{n+1}$ implies, that $a_0=(-1)^na_{n+1}$ and $(-1)^{n+1}a_0<0$. Therefore, there exists a certain index $0\le r\le n$ such that $(-1)^{n+1+s} a_s \le 0$ for $s\le r$ (and is strictly $<0$ for $s=1,\ldots,r-1$), whereas $(-1)^{n+1+q} a_q>0$ for $q>r$.

In particular, the sign of $P_i$ towards $\pm\infty$ can be calculated. Indeed,
\[
\lim_{t\to +\infty}\frac{P_i(t)}{w(t)t^n}=a_i,
\qquad
\lim_{t\to -\infty}\frac{P_i(t)}{w(t)t^n}=(-1)^n a_i.
\]
Hence, for sufficiently large \(T>0\),
\[
\sign P_i(T)=\sign a_i,
\]
and
\[
\sign P_i(-T)=(-1)^n\sign a_i.
\]

Thus, for $i=1,\ldots,n$, an additional zero of $P_i$ in the left outer interval $I_0$ occurs precisely when the sign of $P_i$ at $-\infty$ is opposite to its sign at $x_0$. Similarly, an additional zero in the right outer interval $I_{n+1}$ occurs precisely when
the sign of $P_i$ at $+\infty$ is opposite to its sign at $x_n$.

This gives the complete zero structure of the functions $P_i$: each interior polynomial $P_i$, $i=1,\ldots,n$, has the $n-1$ zeros forced by the sign changes between consecutive nodes, and it may have one further zero in exactly one of the two outer intervals,
depending on the sign of its leading coefficient $a_i$. The endpoint polynomials $P_0$ and $P_{n+1}$ have exactly one simple zero in each bounded interval $(x_{k-1},x_k)$, $k=1,\ldots,n$.

\bigskip
According to the above we now define the index set $\Ji$ to contain all those indices $j$ for which there is a root of $P_i$ in $I_j$, where $j=0,\ldots,n+1$. Then we will have
\begin{equation}\label{Jistardef}
\Ji=\begin{cases}
 \{0,1,\ldots,n\}\setminus \{i\} & \textrm{if} \quad 0\le i<r
\\  & \textrm{and if} \quad i=r \quad \textrm{and} \quad {a_r} \ne 0
\\ \{1,\ldots,n\}\setminus \{i\} & \textrm{if} \quad i=r \quad \textrm{and} \quad {a_r} = 0
\\ \{1,\ldots,n,n+1\}\setminus \{i\} & \textrm{if} \quad r<i\le n+1
\end{cases}. \quad 
\end{equation}

\subsection{Interlacing of roots of $P_i$}
The interlacing argument is purely based on the sign pattern of the functions $P_i$ at the nodes and on the ECHS property. Hence the proof is the same as in the exponential case, with the only additional
observation that, in the Hermite setting, $P_0=(-1)^nP_{n+1}$.
Consequently, the zero set of $P_0$ coincides with that of $P_{n+1}$, so the inclusion of $P_0$ introduces no additional roots. Therefore, the interlacing analysis reduces to the family $P_1,\ldots,P_{n+1}$, exactly as in the exponential case.

We first show that the interlacing property holds pairwise, that is, for every pair $i\neq j$, the zeros of $P_i$ and $P_j$ strictly interlace. This will subsequently imply the cyclic ordering of the root systems of the whole family $P_1,\ldots,P_{n+1}$.

To deal with zeros of different $P_i$, let $y_k^{(i)}$ denote the unique zero of $P_i$ in $I_k$; we have such a zero precisely for $k \in \Ji$. Further, if $i<j$ we denote $\Kij:=\Ji \cap \Jj \setminus \{i+1,\ldots,j-1\}$.

\begin{lemma}\label{l:zerosinsame} Let $\WW=w\Pn$ with a continuous, positive weight function $w$ on $\RR$. Let $1 \le i < j \le n+1$ be arbitrary. If $k \in \Kij$, then the zeros of $P_i$ and $P_j$ satisfy $y_k^{(j)} < y_k^{(i)}$; and if $i<k<j$, then $y_k^{(i)} < y_k^{(j)}  $.
\end{lemma}
Note that $\{i+1,\ldots,j-1\} \subset \Ji \cap \Jj$ and hence the two index sets in the conditions sum up precisely to $\Ji \cap \Jj$.

\begin{proof}
The proof is identical to the corresponding proof in the exponential case. Indeed, the auxiliary functions
\[
G=(-1)^{i+1}P_i+(-1)^{j+1}P_j,
\qquad
H=(-1)^{i+1}P_i-(-1)^{j+1}P_j
\]
have the same sign structure on the intervals $I_k$. Therefore the argument of de Boor and Pinkus \cite{CBoorPinkus}, in the form used in \cite{Exp}, applies verbatim.
\end{proof}

In the next lemma we denote by $f \prec g$ the relation that the (equal in number) zeroes of $f$ and $g$ are ordered so that they can be paired with the zero of $f$ less than that of $g$ in each pair.
\begin{lemma}\label{l:interlacing} Let $\WW=w\Pn$ with a continuous, positive weight function $w$ on $\RR$. For any $1\le i < j \le n+1$, the zeros of $P_i$ and $P_{j}$ strictly interlace, that is, strictly between any two consecutive (simple) zeros of one of the polynomials, there is exactly one (simple) zero of the other, too.

Moreover, the ordering of the $P_i$ is counter-cyclical: $P_r\prec P_{r-1} \prec \dots\prec P_1 \prec P_{n+1} \prec P_n \prec \dots \prec P_{r+1}$,
unless $a_r=0$, in which case $P_r$ is only a degree $n-1$ polynomial, whose $n-1$ zeros are interlacing with the $n$ zeros of any other $P_i$ with $i\ne r$ (and the ordering of the remaining $P_i$ remains counter-cyclical starting with $P_{r-1}$ in this case).
\end{lemma}

\begin{proof} To infer the interlacing and precedence ordering of the root sequences of our polynomials from the zero locations established in Lemma \ref{l:zerosinsame} is a routine derivation, exactly matching that of the classical case, see e.g., \cite{CBoorPinkus} or, for a more detailed deduction, \cite{SzAMH}. Therefore, we save the reader from the details here. The only additional point is the identity $P_0=(-1)^nP_{n+1}$, which shows that $P_0$ introduces no new zeros and hence does not alter the ordering.
\end{proof}

Finally, we get the ordering of all the roots of all the $P_i$.

\begin{corollary}\label{c:order_roots_total}
Let $\WW=w\Pn$ with a continuous, positive weight function $w$ on $\RR$. The roots $y_k^{(i)}$ ($k=0,\dots,n$, $i=1,\dots,n+1$) have the following order:
\begin{align}\label{order_total}
(y_0^{(r)}<)~y_0^{(r-1)}&<\dots<y_0^{(1)}< x_0 <y_1^{(n+1)}<y_1^{(n)}<\dots<y_1^{(2)} \notag
\\ & < x_1 < \dots < x_n <y_{n+1}^{(n)}<\dots<y_n^{(r+1)},
\end{align}
where $y_n^{(r)}$ exists, if and only if $\deg P_r=n$. Moreover, because $P_0=(-1)^nP_{n+1}$, the zeros of $P_0$ coincide with the zeros of $P_{n+1}$.
\end{corollary}

\subsection{Interlacing of roots of the derivatives $P_i'$}\label{sec:Markov}

Now, we are in a position to describe the location of the zeroes of the \emph{derivatives} $P_i'$ of the functions $P_i$.
The previously obtained properties of the zeros of the functions $P_i$ also provide information about the location of the zeros of their derivatives.

In the classical case, namely for algebraic polynomials, V.~A.~Markov proved a well-known interlacing-type lemma concerning the relation between the zeros of polynomials and the zeros of their derivatives, (see \cite{Markov}). It states that if two oscillating polynomials (that is, polynomials having the maximal number of simple zeros) have interlacing zeros, then the zeros of their derivatives interlace as well.

We shall need a similar statement, namely a Markov-type lemma, for Hermite weighted polynomials.

It is important to emphasize that we will need this result not only for two oscillating weighted polynomials, but also in the situation arising from the existence of $P_r$, which may have one less zero in the case $a_r=0$. Thus, we shall also require a version of the statement for one oscillating weighted polynomial and another weighted polynomial having one less zero.

Both versions of this result are known from the works of Milev and Naidenov, see \cite{MN} and \cite{Milev-Naidenov-JAT-2006}, respectively.

\begin{theorem}[Milev-Naidenov]\label{Milev_Naidenov}
Let $\WW_{n}:=\{wp~:~ p\in \Pn\}$, with $w(t)=\exp(-t^2)$, be the space of Hermite weighted polynomials. Let $f$ and $g$ be two oscillating polynomials of $\WW_n$ having zeroes $u_1<\ldots<u_{n}$ and $v_1<\ldots<v_{n}$, respectively and satisfying the interlacing property $u_1\le v_1\le u_2\le\dots\le v_{n-2} \le u_{n}\le v_{n}$.

Then we also have the similar interlacing property for the roots of the derivatives $f', g'$. Moreover, if some of the majorization inequalities $u_i\le v_i$ holds with strict inequality, i.e., $u_i<v_i$, then also the inequalities between the roots of derivatives are strict.
\end{theorem}

\begin{theorem}[Milev-Naidenov]\label{Markov_exp_Pr}
Let $\WW_{n}$ denote the Hermite weighted polynomial system.
Assume that $f$ is an oscillating polynomial of $\WW_n$, and $g\in\WW_{n-1}$, which has $n-1$ real roots. Suppose that the zeroes $u_i$ of $f$ and $v_j$ of $g$ are interlacing, i.e., $u_1\le v_1\le u_2\le\dots \le u_{n-1}\le v_{n-1}\le u_{n}$. Then, we also have a similar, but always strict interlacing property for the roots of $f'$ and $g'$.
\end{theorem}

In the following let us denote the root sequences, arranged in increasing order, of the derivatives $P_i'$ as $W_i$. Recall that $\# W_i=n$ except possibly for $i=r$, where $\# W_r=n-1$ may occur in case $a_r=0$ and $\deg P_r=n-1$. The main goal of our analysis thus far was to obtain the following result.

\begin{corollary}\label{cor:derivativeroots} If $w(t)=\exp(-t^2)$ is the Hermite weight, then the root sequences $W_i$ $(i=1,\ldots,n+1)$ of the derivatives $P_i'$ of the weighted polynomials $P_i$ defined by $P_i:=L_{n+1}|_{I_i}$ on $I_i$ strictly interlace.

If $\deg P_r=n$, then the interlacing is in the order $$W_r\prec W_{r-1}\prec\dots\prec W_1 \prec W_{n+1}=W_0 \prec W_n \prec \dots \prec W_{r+1}.$$

If $\deg P_r=n-1$, then the ordering remains the same except for $W_r$, which has $\# W_r=n-1$ roots, and these zeroes are strictly interlacing with the $n$ zeros of $W_i$ for any $i \ne r$.
\end{corollary}

\section{Maximum place and maximum value on $I_i$}\label{sec:maximum}

The results established in the next two sections are the Hermite-weighted counterparts of the corresponding results proved for the exponentially weighted case in Sections~4 and~5 of \cite{Exp}. In fact, the arguments developed there apply in a much more general setting, namely to weighted polynomial systems $w\Pn$ generated by a twice continuously differentiable, log-concave weight and forming an ECHS. Since the Hermite weight $w(t)=e^{-t^2}$
satisfies these assumptions, most of the proofs carry over verbatim.

The main difference is that, in the present Hermite-weighted setting, there are two outer intervals, $I_0=(-\infty,x_0]$, and $I_{n+1}=[x_n,\infty)$, instead of the single outer interval occurring in the half-line setting of \cite{Exp}. Consequently, besides the weighted polynomial $P_{n+1}$, one also has to consider the polynomial $P_0$, together with the corresponding maxima $m_0$ and maximisers $z_0$.

Whenever the proofs are identical to those in \cite{Exp}, we shall therefore only indicate the necessary modifications related to the additional left outer interval.

Firstly, we find the maximum point of the \emph{weighted} Lebesgue function $L_n$ on each of the intervals $I_i$. That means to find the maximum point of $P_i(\xx,t)=w(t)p_i(t)$ on $I_i$. For sure we know, that $P_i(\xx, t)>0$ on $I_i$, for every $i=0,1,\ldots n+1$. Moreover, $P_i(\xx, x_{i-1})=P_i(\xx, x_{i})=1$ ($i=1,\ldots,n$), and $P_0(x_0)=1$ $P_{n+1}(x_n)=1$. Looking for the maximum place and value of $P_i$ on $I_i$ is equivalent to looking for the maximum place and value of
\[
F_i(\xx,t):=\log |P_i(\xx,t)| = \log |w(t)p_i(t)|,
\]
on $I_i=I_i(\xx)$, $i=0,\ldots,n+1$.

\subsection{Strong locally uniform concavity.}

The following statement is precisely Lemma~7 of \cite{Exp}. Since its proof relies only on the twice continuously differentiability and log-concavity of the weight, it applies verbatim to the present Hermite-weighted setting.

\begin{lemma}\label{l:strict_concave_place_zi}
Let $\WW_n=w\Pn$ be a weighted polynomial system with a twice continuously differentiable and log-concave weight $w$.

Let $K$ be any compact set $K \subset \overline{S}\times \RR$ with the property that it avoids zeros of $P_i$ $(i=0,\ldots,n+1)$, that is, $P_i(\xx,t) \ne 0$ for $(\xx,t)\in K$. Then $F_i$ is strictly concave on $K$ in the variable $t$, and there exists a constant $c=c(K)>0$ such that $(F_i)''_{tt}(\xx,t)<-c$ uniformly for $(\xx,t)\in K$.
\end{lemma}

So at this point we can define the interval maxima $z_i(\xx)$ as a unique maximum point of $F_i$ or, equivalently, of $P_i$, on $I_i$. Given that $P_i(\xx, x_{i-1})=P_i(\xx, x_{i})(=1)$ for all $i=1,\ldots,n$, strict log-concavity immediately implies that $z_i \in \intt I_i$ whenever $1\le i \le n$. Therefore we obtain the following.

\begin{corollary}\label{c:ziininttIi}
If the weight $w$ is twice continuously differentiable and log-concave, then $z_i(\xx) \in \intt I_i$, and hence $m_i(\xx)>1$ for every $i=1,\dots,n$.
\end{corollary}

\subsection{The critical points $w_i(\xx)$.}

As in the exponential case \cite{Exp}, we shall work with the maximal intervals $J_i$ containing $I_i$ on which the weighted polynomial $P_i$ does not change sign. More precisely, $J_i$ is the maximal interval containing $I_i$ such that $P_i(\xx,t)>0$ for all $t\in \intt J_i$, or equivalently, $F_i(\xx,t)>-\infty$ for all $t\in \intt J_i$, while $P_i$ vanishes at the endpoints of $J_i$.

For $i=2,\ldots,n$, the intervals $J_i$ are defined exactly as in the exponential setting, namely $J_i:=(y_{i-1}^{(i)},y_{i+1}^{(i)})$. The interval $J_{n+1}$ is also unchanged and it is determined by $J_{n+1}:=(y_{n}^{(n+1)},\infty)$. The only additional interval in the Hermite-weighted setting is $J_0:=(-\infty,y_1^{(0)})=(-\infty,y_{n+1}^{(0)})$.

Moreover, the description of $J_1$ becomes simpler than in \cite{Exp}. Indeed, since Property~(P) holds with $\de_0=\de_n=1$, the polynomial $P_1$ always possesses a zero to the left of $I_1$ and another one to its right. Consequently, $J_1=(y_0^{(1)},y_2^{(1)})$ without any further case distinction.

So we have defined the maximal open intervals $J_i=J_i(\xx)$ for any node systems $\xx \in S$ such that by construction $[x_{i-1},x_{i}]=I_i\subset J_i$, $P_i>0$ on $J_i$, and the limits of $P_i(\xx,\cdot)$ towards the endpoints of $J_i$ are zero. In other words, $F_i=\log P_i$ tends to $-\infty$ towards the endpoints of $J_i$, whilst $F_i$ is strictly concave on $J_i$. Therefore $F_i(\xx,\cdot)$  has a unique strict maximum point $w_i=w_i(\xx)$ inside $\intt J_i$. Moreover, $F_i(\xx,\cdot)$ is differentiable, and the maximum point must satisfy $(F_i)'_t(\xx,w_i(\xx))=0$. In view of strict concavity, we can identify $w_i(\xx)$ as the unique critical point of $F_i$ on $J_i$. When $z_i(\xx) \in \intt I_i(\xx)$ -- in particular, for the Hermite weight certainly for $1 \le i \le n$ -- the interval maxima $z_i(\xx)$ and the critical points $w_i(\xx)$ coincide. However, for the outer intervals it is well possible that $z_0=x_0<w_0<y_1^{(0)}=y_1^{(n+1)}$, or $y_n^{(n+1)} < w_{n+1} < x_n =z_{n+1}$, and $z_0$ or $z_{n+1}$ are \emph{not} a critical point. We only know that $w_0(\xx)\in J_0(\xx)=(-\infty, y_1^{(0)})$ and $w_{n+1}(\xx) \in J_{n+1}(\xx) = (y_n^{(n+1)},\infty)$.

The next step is to obtain a local uniform localization of the critical points $w_i(\xx)$. Since the argument relies on the strong local uniform concavity of the functions $F_i$ and not on any specific feature of the exponential weight, the proof is identical to that of the corresponding result in \cite{Exp}. We therefore omit the proof and only state the conclusion.

\begin{lemma}\label{l:localization-wi} Let $w$ be a twice continuously differentiable log-concave weight function.

If $\ba \in S$ and $0\le i \le n+1$ is fixed, then there exists $\de>0$ and a compact interval $N=[\alpha, \beta] \subset J_i(\ba)$, such that for $\xx \in B(\ba,\de)$ we also have $w_i(\xx) \in N$.

Moreover, $w_i(\xx)$ is continuously differentiable, and its derivative is uniformly bounded in $B(\ba,\de)$.
\end{lemma}

The assertion is essentially trivial, once we observe that $w_i(\xx)$ is the unique solution of the implicit functional equation $(F_i)'_t(\xx,w_i(t))=0$, and establish that the derivative with respect to $t$, that is, $(F_i)''_{tt}$, is continuous and separated from 0. The below detailed calculus is provided only for certain technical advantages in the further calculations.

\subsection{Lipschitz continuity of the interval maxima $z_i(\xx)$.}

The proof of the Lipschitz continuity of the maximisers $z_i(\xx)$, $i=0,1,\ldots,n+1$, is identical to the proof given in \cite{Exp}. The only new element is the additional endpoint maximiser $z_0$, which arises from the left outer interval $I_0=(-\infty,x_0]$. Since the original argument is brief and the
necessary modifications are minimal, we include the proof for completeness.

\begin{lemma}[Lipschitz continuity of the maximum location on $I_i$] Let $w$ be a twice continuously differentiable log-concave weight function and assume that $w\Pn$ forms an ECHS.

Then the Lebesgue function $L_n(w,\xx,t)$ has a unique maximum point $z_i(\xx)$ in each intervals $I_i$ for $i=0, 1,\ldots,n+1$. This maximum is a strict maximum for $I_i$ and depends on $\xx$ in a locally Lipschitz continuous manner.
\end{lemma}

\begin{proof} That $z_i(\xx)$ exists and is unique, and moreover it is a strict maximum point for $I_i$, follows directly from the strict concavity of $F_i:=\log P_i$. However, we do not always have that $z_i(\xx)=w_i(\xx)$, as $w_i(\xx)\in N \subset J_i$ may be off $I_i$. In our setup, this may happen for $i=0$ and $i=n+1$, when possibly $w_0>z_0=x_0$ or $w_{n+1}<z_{n+1}=x_n$.
Nevertheless, it is still obvious that
$z_0(\xx)=\min(x_0,w_0(\xx))$,
respectively, $z_{n+1}(\xx)=\max(x_n, w_{n+1}(\xx))$.

So we recall that if $g, f : K \to \RR$ are Lipschitz continuous functions on a domain $K$, then $h:=\max(f,g)$ is Lipschitz continuous, too.
\end{proof}

\section{Derivatives of the local maxima functions}

In this section, similarly to the exponentially case, we present the a fundamental consequence: all the $M_i(\xx)=F_i(\xx,z_i(\xx))$ are continuously differentiable. There is no essentially new argument required. The proofs carry over verbatim from the exponential setting \cite{Exp}. The only additional consideration concerns the partial derivatives of the functions $M_0$ and $m_0$, for which one has to apply the same argument symmetrically to the left outer interval $I_0=(-\infty,x_0]$ as was used in \cite{Exp} for the functions $M_{n+1}$ and $m_{n+1}$ on the right outer interval $I_{n+1}=[x_n,\infty)$.

\begin{proposition}\label{prop:Miderivative}
Let $w$ be a twice continuously differentiable log-concave weight function on $\RR$.

At any $\ba \in S$ the interval maxima functions $M_i$ are continuously differentiable for all $i=1,\ldots,n$, and their $j$th partial derivatives are
\begin{equation}\label{Mipartialderivative}
\frac{\partial M_i}{\partial x_j} (\ba)= \frac{\partial F_i }{\partial x_j} (\ba,z_i(\ba)),  \qquad(i=1,\ldots,n,~ j=0,1,\ldots,n),
\end{equation}
for $i=0$ and arbitrary $j=0,1,\ldots,n$
\begin{equation}\label{Mnplus1partialderivative_0}
\frac{\partial M_{0}}{\partial x_j} (\ba) =
\begin{cases}
\dfrac{\partial F_{0}}{\partial x_j}(\ba,z_{0}(\ba)), \qquad \qquad \qquad \textrm{if} \quad \ba \in S\setminus X_-, \quad \text{i.e., } z_{0}(\ba)<a_0
\\[2mm]
\dfrac{\partial F_{0}}{\partial x_j}(\ba,a_0)+(F_{0})'_t(\ba,a_0),  \quad \textrm{if} \quad \ba\in X_-, \quad \text{i.e., }  z_{0}(\ba)=a_0
\end{cases},
\end{equation}
and finally, for $i=n+1$ and arbitrary $j=0,1,\ldots,n$
\begin{equation}\label{Mnplus1partialderivative}
\frac{\partial M_{n+1}}{\partial x_j} (\ba) =
\begin{cases}
\dfrac{\partial F_{n+1}}{\partial x_j}(\ba,z_{n+1}(\ba)), \qquad \qquad \qquad \textrm{if} \quad \ba \in S\setminus X_+, \quad \text{i.e., }  z_{n+1}(\ba)>a_n
\\[2mm]
\dfrac{\partial F_{n+1}}{\partial x_j}(\ba,a_n)+(F_{n+1})'_t(\ba,a_n),  \quad \textrm{if} \quad \ba\in X_+, \quad  \text{i.e., } z_{n+1}(\ba)=a_n
\end{cases}.
\end{equation}
In other words,
\begin{equation}\label{Migradient}
\nabla_{\xx} M_i(\ba) = \nabla_\xx F_i(\ba,z_i(\ba)), \qquad \qquad  \textrm{if} \quad i=1,\ldots,n ,
\end{equation}
\begin{equation}\label{Mnplus1gradient_0}
\nabla_{\xx} M_{0}(\ba) =
\begin{cases}
\nabla_\xx F_{0}(\ba,z_{0}(\ba)), \qquad \qquad \textrm{if} \quad  \ba \in S\setminus X_-, \quad \text{i.e., } z_{0}(\ba)<a_0
\\ \nabla_\xx F_{0}(\ba, a_0)+(F_{0})'_t(\ba, a_0) \ee_1, \quad \textrm{if} \quad \ba\in X_-, \quad \text{i.e., }  z_{0}(\ba)=a_0
\end{cases}
\end{equation}
and
\begin{equation}\label{Mnplus1gradient}
\nabla_{\xx} M_{n+1}(\ba) =
\begin{cases}
\nabla_\xx F_{n+1}(\ba,z_{n+1}(\ba)), \qquad \qquad \textrm{if} \quad \ba \in S\setminus X_+, \quad \text{i.e., }  z_{n+1}(\ba)>a_n
\\ \nabla_\xx F_{n+1}(\ba, a_n)+(F_{n+1})'_t(\ba, a_n) \ee_{n+1}, \quad \textrm{if} \quad \ba\in X_+,  \quad \text{i.e., }  z_{n+1}(\ba)=a_n.
\end{cases}
\end{equation}
\end{proposition}

Note that the above formulae describe a continuous derivative, since whenever $z_{n+1}(\xx) \in \intt I_{n+1}$, then we have $(F_{n+1})'_t(\xx,z_{n+1}(\xx))=0$, hence when $\xx \to \ba$ with $z_{n+1}(\xx)\in \intt I_i(\xx)$ but $z_{n+1}(\xx)\to a_n$, i.e., $z_{n+1}(\ba)=a_n$, then also the endpoint value of $(F_{n+1})'_t(\ba,z_{n+1}(\ba))$ is equal to $\lim_{\xx \to \ba} (F_{n+1})'_t(\xx,z_{n+1}(\xx))=0$ by continuity.

\begin{proposition}\label{prop:Mnplus1derivative}
Let $w$ be a twice continuously differentiable log-concave weight function.
If for an $\ba \in S$ $z_{0}(\ba)=a_0$ holds, i.e., $\ba\in X_-$, then we also have
$$
\dfrac{\partial M_{0}}{\partial x_j} (\ba) =0, \qquad (j=0,\ldots,n).
$$
That is, we also have
\begin{equation}\label{Mnplus1gradientzero_0}
\nabla_{\xx} M_{0}(\ba) = {\bf 0}, \qquad \textrm{if} \quad \ba\in X_-, \quad \text{i.e., }  z_{0}(\ba)=a_0.
\end{equation}

Similarly, if for an $\ba \in S$ $z_{n+1}(\ba)=a_n$ holds, i.e., $\ba\in X_+$, then we also have
$$
\dfrac{\partial M_{n+1}}{\partial x_j} (\ba) =0, \qquad (j=0, 1,\ldots,n).
$$
That is, we also have
\begin{equation}\label{Mnplus1gradientzero}
\nabla_{\xx} M_{n+1}(\ba) = {\bf 0}, \qquad \textrm{if} \quad \ba\in X_+, \quad \text{i.e., }  z_{n+1}(\ba)=a_n.
\end{equation}
\end{proposition}

The proof of the formula for $\frac{\partial  M_{n+1}}{\partial x_j} (\ba)$ under the assumption $z_{n+1}(\ba)=a_n$ is identical to that of Proposition~4 in \cite{Exp}, where the set denoted there by $X^c$ is replaced by $X_+$. Proposition~\ref{prop:fat-hermite} establishes the required fatness property for both $X_+$ and $X_-$. The corresponding formula for  $\frac{\partial  M_{0}}{\partial x_j} (\ba)$, under the assumption $z_0(\ba)=a_0$ then follows by symmetry.

%
%

Now, we can compute the derivatives of the local maxima functions of $P_i(\xx,\cdot)$ on the interval $I_i(\xx)$. Recall $m_i(\xx)=\exp(M_i(\xx))$, and hence $m_i(\xx)=\exp(F_i(\xx,z_i(\xx)))=P_i(\xx,z_i(\xx))=\max_{I_i(\xx)}P_i(\xx,\cdot)$. Then, from the logarithmic derivative formulae in Proposition~\ref{prop:Miderivative} and Proposition~ \ref{prop:Mnplus1derivative} we also have the following.

\begin{proposition}\label{prop:miderivative}
Let $w$ be a twice continuously differentiable log-concave weight function.

At any $\ba \in S$ the interval maxima functions $m_i$ are continuously differentiable for all $i=0,1,\ldots,n+1$, and their $j$th partial derivatives are
\begin{equation}\label{mipartialderivative}
\frac{\partial m_i}{\partial x_j} (\ba)= \frac{\partial P_i}{\partial x_j} (\ba,z_i(\ba)),  \qquad(i=1,\ldots,n,~ j=0,1,\ldots,n),
\end{equation}
for $i=0$ and arbitrary $j=0,1,\ldots,n$
\begin{equation}\label{mnplus1partialderivative_0}
\frac{\partial m_{0}}{\partial x_j} (\ba) =
\begin{cases}
\frac{\partial P_{0}}{\partial x_j}(\ba,z_{0}(\ba)), \qquad \qquad \qquad \textrm{if} \quad \ba\in S\setminus X_-, \quad \text{i.e., } z_{0}(\ba)<a_0
\\[2mm]
\frac{\partial P_{0}}{\partial x_j} (\ba,a_0) +(P_{n+1})'_t(\ba,a_n)=0, \quad \textrm{if} \quad \ba\in  X_-, \quad \text{i.e., } z_{0}(\ba)=a_0
\end{cases},
\end{equation}
and for $i=n+1$ and arbitrary $j=0,1,\ldots,n$
\begin{equation}\label{mnplus1partialderivative}
\frac{\partial m_{n+1}}{\partial x_j} (\ba) =
\begin{cases}
\frac{\partial P_{n+1}}{\partial x_j}(\ba,z_{n+1}(\ba)), \qquad \qquad \qquad \textrm{if} \quad \ba\in S\setminus X_+, \quad \text{i.e., } z_{n+1}(\ba)>a_n
\\[2mm]
\frac{\partial P_{n+1}}{\partial x_j} (\ba,a_n) +(P_{n+1})'_t(\ba,a_n)=0, \quad \textrm{if} \quad \ba\in  X_+, \quad \text{i.e., }  z_{n+1}(\ba)=a_n.
\end{cases}
\end{equation}
In other words,
\begin{equation}\label{migradient}
\nabla_{\xx} m_i(\ba) = \nabla_\xx P_i(\ba,z_i(\ba)), \qquad \qquad  \textrm{if} \quad i=1,\ldots,n ,
\end{equation}
\begin{equation}\label{mnplus1gradient_0}
\nabla_{\xx} m_{0}(\ba) =
\begin{cases}
\nabla_\xx P_{0}(\ba,a_{0}), \qquad \qquad \textrm{if} \quad \ba\in S\setminus X_-, \quad \text{i.e., } z_{0}(\ba)<a_0
\\[2mm] \nabla_\xx P_{0}(\ba,a_0)+(P_{0})'_t(\ba,a_0) \ee_1={\bf 0}, \quad \textrm{if} \quad \ba\in X_-, \quad \text{i.e., } z_{0}(\ba)=a_0,
\end{cases}
\end{equation}
and
\begin{equation}\label{mnplus1gradient}
\nabla_{\xx} m_{n+1}(\ba) =
\begin{cases}
\nabla_\xx P_{n+1}(\ba,a_{n}), \quad \qquad \textrm{if} \quad \ba\in S\setminus X_+, \quad \text{i.e.,} \quad z_{n+1}(\ba)>a_n
\\[2mm] \nabla_\xx P_{n+1}(\ba,a_n)+(P_{n+1})'_t(\ba,a_n) \ee_{n+1}={\bf 0}, \quad \textrm{if} \quad  \ba\in X_+, \quad \text{i.e., } z_{n+1}(\ba)=a_n.
\end{cases}
\end{equation}
\end{proposition}

\section{Consequences for the signs of partial derivatives}

The next result is mentioned without proof in \cite{K-AMH}, and it is called only ``an observation''. According to our study, it can be traced back to the work of Kilgore and Cheney \cite{Kilgore-Cheney}. However, Kilgore and Cheney use that ${\bf 1}$ belongs to the given ECHS, which entails\footnote{Once $z_i$ can be defined uniquely, which does not follow from these conditions alone, but require also that the derivative system be ECHS.} $z_i \in \intt I_i$, for all $i$, and therefore they do not have the respective formula for the cases when $z_0(\xx)=x_0$, i.e., $\xx\in X_-$ or $z_{n+1}=x_n$, i.e., $\xx\in X_+$.

The following statement is the Hermite-weighted analogue of Corollary~5 in \cite{Exp}. Apart from the inclusion of the corresponding formulas for the partial derivatives of $P_0$, no changes are required. Since the proof is identical to that of \cite{Exp}, we merely state the result.

\begin{lemma}[Magic formula]\label{l:magic}
Let $w$ be a twice continuously differentiable log-concave weight function. Let $\xx \in S$. Then, for every $0\le j\le n$ and $1\le i\le n$ we have
\begin{equation}\label{magic-middleform}
\frac{\partial P_i}{\partial x_j}(\xx,t)= - h_j(\xx,t) (P_i)'_t(\xx,x_j), \qquad (t \in \RR).
\end{equation}
Moreover, for every $0\le j\le n$ and $1\le i\le n$ it holds
\begin{equation}\label{der_identity}
\frac{\partial m_i}{\partial x_j}(\xx)=-h_j(\xx,z_i(\xx))(P_i)'_t(\xx,x_j),
\end{equation}
whereas for the derivative of $m_{0}$, we have for every $0\le j\le n$ that
\begin{equation}\label{der_final_0}
\frac{\partial m_{0}}{\partial x_j}(\xx)=\begin{cases}
-h_j(\xx,z_{0}(\xx))(P_{0})'_t(\xx,x_j), & \text{ if } \quad \xx \in S\setminus X_-, \quad \text{ i.e.,  }  z_0(\xx)<x_0
\\ 0,  & \text{ if } \quad \xx \in X_-, \quad \text{ i.e., }  z_0(\xx)=x_0,
\end{cases}
\end{equation}
and for the derivative of $m_{n+1}$, we have for every $0\le j\le n$ that
\begin{equation*}
\frac{\partial m_{n+1}}{\partial x_j}(\xx)=\begin{cases}
-h_j(\xx,z_{n+1}(\xx))(P_{n+1})'_t(\xx,x_j), & \text{ if } \quad \xx \in S\setminus X_+, \quad \text{ i.e., } z_{n+1}(\xx)>x_n
\\ 0,  & \text{ if } \quad \xx \in X_+, \quad \text{ i.e., }  z_{n+1}(\xx)=x_n .
\end{cases}
\end{equation*}
\end{lemma}

\begin{corollary}\label{c:countermonotonicity}
Let $w$ be a twice continuously differentiable log-concave weight function.

\begin{itemize}
\item
If $\xx\in X$, then
\begin{equation}\label{signchange}
\frac{\partial m_i}{\partial x_i}(\xx)>0  \qquad \text{ and } \qquad
\frac{\partial m_{i+1}}{\partial x_i}(\xx)<0, \qquad (i=0,\dots,n).
\end{equation}

\item
If $\xx \in X_+\setminus X_-$, then \eqref{signchange} still holds for every $i<n$, but for $i=n$ and for any $\xx\in X_+$ we have
\begin{equation}\label{signchangeXc}
\frac{\partial m_n}{\partial x_n}(\xx)>0  \qquad \text{ and } \qquad
\frac{\partial m_{n+1}}{\partial x_n}(\xx)=0.
\end{equation}
\item
If $\xx \in X_-\setminus X_+$, then \eqref{signchange} still holds for every $i>0$, but for $i=0$ and for any $\xx\in X_-$ we have
\begin{equation}\label{signchangeXc_2}
\frac{\partial m_{0}}{\partial x_0}(\xx)=0 \qquad \text{ and } \qquad \frac{\partial m_1}{\partial x_0}(\xx)<0.
\end{equation}
\end{itemize}
In particular, we get that for any $\xx\in X\cup (X_+\setminus X_-)=S\setminus X_-$
\begin{equation}\label{signchangeX+}
\frac{\partial m_{0}}{\partial x_0}(\xx)>0,
\end{equation}
and for any $\xx\in X\cup (X_-\setminus X_+)=S\setminus X_+$
\begin{equation}\label{signchangeX-}
\frac{\partial m_{n+1}}{\partial x_n}(\xx)<0.
\end{equation}

\end{corollary}

\begin{proof}
The properties listed in \eqref{signchange} and \eqref{signchangeXc} are established exactly as in Corollary~6 of \cite{Exp}. The remaining sign properties in \eqref{signchangeXc_2} follow by symmetry from those in \eqref{signchangeXc}.
\end{proof}

Observe that in particular we obtain the following corollary, which will gain importance in Section \ref{sec:unicity}.

\begin{corollary}\label{c:n-n-monotonicity} Let $w$ be a twice continuously differentiable log-concave weight function.

Then $m_{n+1}(\xx)$ is monotonically non-increasing with $x_0<x_1<\dots<x_{n-1}$ fixed and $x_n$ changing in $(x_{n-1},\infty)$.

Similarly, $m_{0}(\xx)$ is monotonically non-decreasing with $x_1<\dots<x_{n-1}<x_n$ fixed and $x_0$ changing in $(-\infty,x_1)$.
\end{corollary}

So far, we have established only those properties of the partial derivatives that are identical to, or direct analogues of, the corresponding results in the exponential setting \cite{Exp}. The key difference is that, in the exponential case, $X^c=\{\xx \in S: m_{n+1}(\xx)=1\}$, whereas in the Hermite-weighted setting $X^c=X_-\cup X_+$.

Consequently, the boundary structure of the regular set is more intricate, and the proof of our main results requires several additional inequalities for the partial derivatives. These are established in the following lemma.

\begin{lemma}\label{l:dmidxn_2}
Let $w$ be a twice continuously differentiable log-concave weight function on $\RR$ and let $\xx=(x_0,x_1,\dots, x_n)\in S$.

If $\xx\in X$, then 
\begin{equation}\label{signchange_outer_m}
\frac{\partial m_{0}}{\partial x_n}(\xx)< 0,
 \qquad \text{ and } \qquad
\frac{\partial m_{n+1}}{\partial x_0}(\xx)>0.
\end{equation}

If $\xx \in X_+$, then 
\begin{equation}\label{signchange_mi}
\frac{\partial m_i}{\partial x_n}(\xx)>0, \quad \text{ for } \quad  i=1,\dots, n,
\qquad \text{ and } \qquad
\frac{\partial m_{0}}{\partial x_n}(\xx)\ge 0.
\end{equation}

Finally, if $\xx\in X_-$, then 
\begin{equation}\label{signchange_mi_X-}
\frac{\partial m_i}{\partial x_0}(\xx)<0, \quad \text{ for } \quad  i=1,\dots, n,
\qquad \text{ and } \qquad
\frac{\partial m_{n+1}}{\partial x_0}(\xx)\le 0.
\end{equation}
\end{lemma}

\begin{proof}
First, several statements are symmetric\footnote{If for a node system $\xx$, we get that $z_{n+1}(\xx)<x_n$, respectively $z_{n+1}(\xx)=x_n$, then for $-\xx$ we get that $z_0(-\xx)<-x_n$, resp. $z_0(-\xx)=-x_n$.}. More precisely, the two inequalities in \eqref{signchange_outer_m} are symmetric and the properties listed in \eqref{signchange_mi} and those in  \eqref{signchange_mi_X-} are symmetric, too.

It means that it is enough to prove the first inequality of \eqref{signchange_outer_m}, and \eqref{signchange_mi}.

Concerning the first inequality of \eqref{signchange_mi}, let $\xx$ be an element of $X_+$, i.e., $z_{n+1}(\xx)=x_n$, ($m_{n+1}(\xx)=1$). Then $P_{n+1}(\xx,x_n)=1$, and $P_{n+1}(\xx,t)<1$ for $t>x_n$. It follows that $(P_{n+1})'_t(x_n)\le 0$. Consider any other $P_i$ with $0<i\le n$. Then for $t>x_{n}$, $|P_i(\xx,t)|<P_{n+1}(\xx,t)$. Consider derivatives close to $x_n$, where we have also $|P_i|\approx 1$, hence $|P_i|$ differentiable. Note that $(|\cdot|')=\sign$ and hence $|P_i|'_t(x_n)=\sign(P_i(x_n)) (P_i)'_t(x_n)$. We therefore obtain
\begin{align*}
|P_i|'_t(x_n) & -(P_{n+1})'_t(x_n)  = \sum_{j=0}^n \left(\sign(P_i(\xx,x_n)) \vji - \ve_{j,n+1}\right)h_j'(x_n)
\\ & = \sum_{j=0}^n \left((-1)^{n-i} (-1)^{i+j+1+\chi_{i\le j}} - (-1)^{j+n+1+1+\chi_{n+1\le j}}\right)h_j'(x_n)
\\ & = \sum_{j=0}^n \left((-1)^{\chi_{i\le j}} +1\right)(-1)^{n+j+1} h_j'(x_n)=2 \sum_{j=0}^{i-1} (-1)^{n+j+1} h_j'(x_n).
\end{align*}
Now, $h_j$ has a maximal number of $n$ zeroes at the nodes $x_i$ with $i\ne j$, whence these nodes are simple, and cannot be zeroes of the derivatives $h_j'$. Further, $h_j$ changes sign from positive to negative at $x_{j+1}$, and so on, from $(-1)^{n+j+1}$ to $(-1)^{n+j}$ at $x_n$, so that $\sign h_j'(x_n)= (-1)^{n+j}$. It follows that all the terms in the above sum are negative, and therefore
$$
|P_i|'_t(x_n)  -(P_{n+1})'_t(x_n) < 0.
$$
Thus we obtained a little stronger statement than the mere observation that $|P_i(\xx,t)|<P_{n+1}(\xx,t)$. Moreover, in this calculation we did not use $\xx \in X_+$. (It is true that when $z_{n+1} \in \intt I_{n+1}$ then $(P_{n+1})'_t(x_n)>0$, and our inequality is not suitable for determining the sign of $|P_i|'_t(x_n)$, but it is still  valid formula.)
Now if $\xx \in X_+$, then $|P_i|'_t(x_n)< (P_{n+1})'_t(x_n) \le 0$ (for $P_{n+1}(x_n)=1$ and then in the interior points of $I_{n+1}$ $P_{n+1}'<0$) so that we obtain $|P_i|'_t(x_n)<0$.

Next we appeal to \eqref{der_identity} of Lemma \ref{l:magic}, according to which
$$
\frac{\partial m_i}{\partial x_n}(\xx) = -h_n(\xx,z_i(\xx))(P_i)'_t(\xx,x_n) \qquad (i=1,\ldots,n).
$$
Note that the formula fails for $i=n+1$ as we had set $\xx \in X_+$, where the partial derivative is just 0. We know that $h_n>0$ on $\intt I_n$, and it changes sign at nodes and has constant sign on all further $I_i$, arriving to $(-1)^{n-i}$ on $I_i \ni z_i$, i.e.,
\begin{equation}\label{hn_sign}
\sign h_n(\xx,z_i(\xx))=(-1)^{n-i}.
\end{equation}

That settles the sign of the first factor. As for the second, $$0>|P_i|'_t(x_n)=\sign (P_i(x_n)) (P_i)'_t(x_n) =(-1)^{n-i} (P_i)'_t(x_n)$$ resulting in $\sign (P_i)'_t(x_n) = (-1)^{n-i+1}$. Combining the two findings, we are led to
$$
\sign\left(\frac{\partial m_i}{\partial x_n}(\xx)\right) =(-1) (-1)^{n-i} (-1)^{n-i+1} =+1,
$$
for every $0<i\le n$.

Now, we consider the first inequality of \eqref{signchange_outer_m} and the second inequality of \eqref{signchange_mi}. Concerning $\partial m_0/\partial x_n(\xx)$, we recall the formula \eqref{der_final_0}, which states that
\[
\frac{\partial m_0}{\partial x_n}(\xx)=-h_n(\xx,z_0(\xx))(P_0)'_t(\xx,x_n)+\chi_{z_0(\xx)=x_n} (P_0)'(\xx,z_0(\xx)).
\]
It is clear that $z_0(\xx)$ differs from $x_n$, therefore the formula simplifies to $$\frac{\partial m_0}{\partial x_n}(\xx)=-h_n(\xx,z_0(\xx))(P_0)'_t(\xx,x_n).$$ Moreover, we know, that $P_0=(-1)^nP_{n+1}$, and hence $(P_0)'_t(\xx,x_n)=(-1)^n(P_{n+1})'_t(\xx,x_n)$.
Firstly, we consider the second inequality of \eqref{signchange_mi}, i.e., we assume that $\xx\in X_+$, i.e., $z_{n+1}(\xx)=x_n$. Then, we know that $(P_{n+1})'_t(\xx,x_n)\le 0$ and the equation \eqref{hn_sign} holds also for $i=0$, i.e., $\sign h_n(\xx,z_0(\xx))=(-1)^n$. If $\xx\in X_+\cap X_-$, i.e., $x_0=z_0(\xx)$, then $h_n(\xx,z_0(\xx))=h_n(\xx,x_0)=0$, which guarantees, that $\frac{\partial m_0}{\partial x_n}(\xx)=0$. Otherwise, i.e., when $\xx\in X_+\setminus X_-$, then we get, that
\[
\sign\left(\frac{\partial m_0}{\partial x_n}(\xx)\right)=(-1)(-1)^n(-1)^n(P_{n+1})'_t(\xx,x_n)\ge 0.
\]

Concerning the first inequality of \eqref{signchange_outer_m}, let $z_{n+1}(\xx)>x_n$. Then, $(P_{n+1})'_t(x_n)> 0$, which implies, that
\[
\sign\left(\frac{\partial m_0}{\partial x_n}(\xx)\right)=(-1)(-1)^n(-1)^n(P_{n+1})'_t(x_n)<0.
\]
\end{proof}

\begin{remark}\label{r:strict_monotonicity}
The proof above yields a slightly stronger statement than the one stated in the second inequalities of \eqref{signchange_mi} and \eqref{signchange_mi_X-}. Namely,
\[
\frac{\partial m_0}{\partial x_n}(\xx)>0 \qquad\text{for every }\xx\in \intt(X_+)\setminus X_-; \qquad \frac{\partial m_{n+1}}{\partial x_0}(\xx)<0 \qquad\text{for every }\xx\in \intt(X_-)\setminus X_+.
\]

Indeed, if $\xx\in \intt(X_+)$, then $x_n>x_n^*=\phi_+(x_0,\ldots,x_{n-1})$, and therefore $(P_{n+1})'_t(\xx,x_n)<0$ (see Lemma \ref{l:xstar-hermite} and Proposition \ref{prop:fat-hermite}). Moreover, if $\xx\notin X_-$, then $z_0(\xx)<x_0$, whence $h_n(\xx,z_0(\xx))\neq 0$. Consequently, the sign computation in the proof becomes strict, yielding $\frac{\partial m_0}{\partial x_n}(\xx)>0$.

Similar argument shows that $\frac{\partial m_{n+1}}{\partial x_0}(\xx)<0$ for every $\xx\in \intt(X_-)\setminus X_+$.
\end{remark}

\section{Linear independence}\label{sec:linearindependence}

The present section follows the classical ideas of Kilgore \cite{Kilgore1985} and de Boor--Pinkus \cite{CBoorPinkus}, and we retain the notation introduced in the corresponding section of \cite{Exp}. Nevertheless, the Hermite-weighted setting is not merely a verbatim repetition of the exponential case. Besides the function $m_{n+1}$, one now has to deal with the additional outer function $m_0$, associated with the left outer interval $I_0$. As a consequence, several arguments require a separate treatment of this second exceptional case.

\begin{definition}\label{Amatrix} For arbitrary $\vv \in S$ we define the column vectors
$$
\ba_i:=\ba_i(\vv):= \nabla_{\xx} m_i(\vv) =\left[\frac{\partial m_i}{\partial x_j}(\vv)\right]_{j=0}^n \qquad (i=0,\ldots,n+1),
$$
and the matrix
$$
A:=A(\vv):=\left[\ba_i(\vv)\right]_{i=0}^{n+1}= \left[\frac{\partial m_i}{\partial x_j}(\vv)\right]_{j=0, i=0}^{n,n+1}.
$$
Further, we denote by $A_i$ the square matrices, obtained  by deleting the $i$th column $\ba_i$ from $A$, and for the (Jacobi) determinants of the $A_i$ we write $D_i:=\det A_i$.
\end{definition}

\begin{definition}\label{Cmatrix} For arbitrary $\vv \in S$
we define the column vectors
$$
\cc_0:=\cc_{0}(\vv):=
\left[\frac{\partial P_{0}}{\partial x_j} (\vv,w_i(\vv)) \right]_{\dg{i=1}}^n, \qquad
\cc_{n+1}:=\cc_{n+1}(\vv):= \left[\frac{\partial P_{n+1}}{\partial x_j} (\vv,w_i(\vv)) \right]_{j=0}^n ,
$$
and the matrix $C$ obtained from $A$ replacing its first and last column by $\cc_0$ and $\cc_{n+1}$ respectively:
$$
C:=C(\vv):=\left[c_{j,i} \right]_{j=0, i=0}^{n,n+1}:=\left[\frac{\partial P_i}{\partial x_j} (\vv,w_i(\vv))\right]_{j=0, i=0}^{n,n+1}=\left[\cc_0,\cc_1, \ldots,\cc_n, \cc_{n+1}\right].
$$
Further, we denote by $C_i$ the square matrices, obtained  by deleting the $i$th column $\cc_i$ from $C$.
\end{definition}

For the interior intervals $I_1,\ldots,I_n$, the argument is identical to that of the exponential case, see \cite{Exp}. Indeed, whenever $1\le i\le n$, we have $z_i(\xx)=w_i(\xx)$, and hence $\frac{\partial P_i}{\partial x_j} (\xx,z_i(\xx))=\frac{\partial P_i}{\partial x_j} (\xx,w_i(\xx))$.
Moreover, the same extends to $i=0$ and $i=n+1$ for $\xx\in X$, and therefore we obtain $A(\xx)=C(\xx)$.

The only additional feature of the Hermite-weighted setting is the presence of the second outer interval $I_0$. As in the exponential case, the column corresponding to $m_{n+1}$ vanishes on $X_+$. By symmetry, the column corresponding to $m_0$ vanishes on $X_-$. Hence, for $\xx\in X^c=X_-\cup X_+$, at least one of the first and last columns of $A$ is zero.

As in \cite{Exp}, we therefore work with the matrices $C_i$ instead of $A_i$.
By \eqref{magic-middleform_0},
\[
c_{j,i}=-h_j(\xx,w_i(\xx))(P_i)'_t(\xx,x_j)=
\frac{w(-w_i)\prod_{l=0}^{n}(w_i-x_l)}{w(-x_j)\prod_{l=0,\,l\neq j}^{n}(x_j-x_l)}
\frac{(P_i)'_t(\xx,x_j)}{x_j-w_i}.
\]

Following the classical reduction of Braess, Kilgore, Cheney and de Boor--Pinkus\footnote{The above reduction was originally suggested by Dietrich Braess; see \cite{Kilgore-Cheney,K-AMH,CBoorPinkus}.}, we divide the $i$th column by the nonzero factor $w(-w_i)\prod_{l=0}^{n}(w_i-x_l)$,
and multiply the $j$th row by $w(-x_j)\prod_{l=0,\,l\neq j}^{n}(x_j-x_l)$. Since these are nonzero row and column scalings, the singularity of the matrices $C_i$ is preserved. Consequently, it is sufficient to study the corresponding submatrices $Q_i$ of
\begin{equation}\label{Qqdef}
Q=[q_i(x_j)]_{j=0,i=0}^{n,n+1},
\qquad
q_i(t):=\pm\frac{(P_i)'_t(\xx,t)}{t-w_i}.
\end{equation}

Although the signs of
\[
\frac{(P_i)'_t(\xx,t)}{t-w_i}\Big|_{t=w_j}
\]
can be computed explicitly, it is more convenient to follow Kilgore and choose the signs in the definition of $q_i$ so that
\[
q_i(w_0)>0,\qquad i=0,\ldots,n+1.
\]

\begin{lemma}\label{l:Kilgore_assump_Hermite}
Consider the Hermite weight on the real line. With the notations above, the functions $q_i$ defined in \eqref{Qqdef} satisfy the Kilgore separation assumptions with respect to the ordered system
\[
w_0'<w_1'<\cdots<w_n'<w'_{n+1}
\]
where with a suitable small $\varepsilon>0$
\[
w_0':=w_0-\varepsilon, \qquad w_i':=w_i \quad (i=1,\ldots,n) \quad w'_{n+1}:=w_{n+1}+\varepsilon
.
\]
More precisely:
\begin{enumerate}
\item[1)] $q_i(w'_j)\ne 0$ for every $i,j=0,\ldots,n+1$;
\item[2)] $q_i$ has a unique simple zero in each interval $[w_j',w'_{j+1}]$ with $j\neq i-1,i$;
\item[3)] $q_i$ has no zero in $[w'_{i-1},w'_{i+1}]$, $i=1,\ldots, n$. Moreover, $q_0$ and $q_{n+1}$ do not have any zero in $(-\infty,w'_1]$, and in $[w'_n,\infty)$, respectively.
\end{enumerate}
\end{lemma}

\begin{proof}
The proof follows almost verbatim the exponential case (see Lemma~11 in \cite{Exp}), except for the ``outer functions'' $q_0, q_{n+1}$, which are derived from $(P_0)'_t=(-1)^n(P_{n+1})'_t$, and hence cause a little technical difficulty. All the properties follow from the cited proof for the very system of points $w_0,w_1,\ldots,w_n,w_{n+1}$, except that we encounter the endpoint vanishing $q_0(w_{n+1})=0$ and $q_{n+1}(w_0)=0$. To avoid this, while keeping all other properties, is possible by a little perturbation decreasing $w_0$ and increasing $w_{n+1}$: this leads to the definition of the perturbed point system of $\{ w_i'\}$. Then the first resp. last intervals between the points $w'_i$ are increased a little, and they will still contain the zeroes they did contain so far; but with small enough perturbation, no new zeroes may be contained. So this perturbed system $w_0',\ldots,w_{n+1}'$ will satisfy even assumption 1).

To see the last assertion -- second sentence in 3) -- recall that
$q_0$ and $q_{n+1}$ are obtained from the same derivative by dividing out two different simple zero factors, namely for $w_0$ and $w_{n+1}$, respectively. Consequently their zero sets are
\[
Z(q_0)=Z((P_{n+1})'_t)\setminus\{w_0\},
\qquad
Z(q_{n+1})=Z((P_{n+1})'_t)\setminus\{w_{n+1}\}.
\]
Moreover, all these zeros are simple. As the zero $w_0$ of $(P_{n+1})'_t$ is removed by the division defining $q_0$, and similarly the zero $w_{n+1}$ is removed in the definition of $q_{n+1}$, we get
\[
q_0(w_0)\ne 0,
\qquad
q_{n+1}(w_{n+1})\ne 0,
\]
while the remaining derivative zeros persist as simple zeros of $q_0$ and $q_{n+1}$, respectively. In particular, $q_0$ has no zero to the left of $w_1$, while $q_{n+1}$ has no zero to the right of $w_n$, (and each has one simple zero in every remaining admissible gap).
In particular, decreasing $w_0$ to $w_0'$ will keep these properties of $q_0$: $q_0(w_0')\ne 0$ and $q_0 \ne 0$ in $(-\infty,w_0']$; and similarly, $q_{n+1}(w_{n+1}') \ne 0$ and $q_{n+1} \ne 0$ in $[w_{n+1},\infty)$ either.
\end{proof}

From this point onward, the argument closely follows that of the exponential case. Since the proofs of the following lemma and corollary are short, we include them here for the reader's convenience.

We note that, once the root distribution described in Lemma~\ref{l:Kilgore_assump_Hermite} has been established, and since all the functions $q_i$ belong to the rank $n+1$ ECHS $\WW_n$, we are in a position to apply a general result of Kilgore.\footnote{The proof follows the classical argument, where the Chebyshev property was immediate from degree considerations for the corresponding polynomials. The elegant proof in the classical setting is due to de Boor and Pinkus \cite{CBoorPinkus}.} More precisely, Proposition~1 of \cite{Kilgore1985} states that, under these assumptions, the system $\{q_0,\ldots,q_{n+1}\}\setminus \{q_k\}$ is linearly independent.

\begin{lemma}[Kilgore]\label{l:Kilgore_orig}
Let $w_0, w_1, \dots, w_{n+1}$ be real numbers with $w_0<w_1<\dots< w_{n+1}$ and $q_0,\dots, q_{n+1}$ be functions which lie in a rank $n+1$ ECHS. Further, assume that $q_0,\dots, q_{n+1}$ satisfy the three conditions 1)-3) listed in Lemma \ref{l:Kilgore_assump_Hermite}.
Then, for any $k\in \{0,\dots, n+1\}$, the set $\{q_0,\dots, q_{n+1}\}\setminus \{q_k\}$ forms a linearly independent set of functions.
\end{lemma}

\begin{lemma}\label{lem:keylemma} Let $q_i$ ($i=0,\ldots,n+1$) belong to an order $n+1$ ECHS and satisfy the three assumptions on their roots as listed in Lemma \ref{l:Kilgore_assump_Hermite}.

Assume that $\sum_{i=0}^{n+1} \ai q_i(x_j)=0$ for all $j=0,\ldots,n$ such that $\al_k=0$ for some $k$ $(0 \le k \le n+1)$. Then we must have that $\al\equiv 0$.
\end{lemma}

\begin{proof}
We already know, that $q_0,\dots, q_{n+1}$ belong to the $n+1$ dimensional ECH system $\WW_{n}=\{e^{-t^2}p(t) : p\in\mathcal{P}_{n}\}$. Since $\sum_{i=0, i\ne k}^{n+1} \ai q_i(x_j)=0$ for all $j=0,\ldots,n$, it means, that there are $n+1$ different zeros of this function belonging to $\WW_{n}$. Therefore, $\sum_{i=0,i\not= k}^{n+1} \ai q_i\equiv 0$.

From here, linear independence of the system $\{q_i\}_{i=0, i\ne k}^{n+1}$, given by Lemma \ref{l:Kilgore_orig}, guarantees that the only possibility is $\al \equiv 0$.
\end{proof}

Consequently, dropping any column from $C$ (or equivalently from $Q$), the remaining columns cannot have a zero linear combination with any other coefficients than with identically zero $\ai$ ($i=0,\ldots,n+1, i\ne k$). That is, $Q_k$ are nonsingular matrices.

That means for $\xx\in X$ the nonsingularity of the matrices $A_k$, as for $\xx \in X$ we have $A=C$.
Moreover, if $\xx\in X_+$, then dropping the last column from $A$, (similarly, if $\xx\in X_-$, then dropping the first column of $A$) we still get the same conclusion, as $A_{n+1}=C_{n+1}$, (respectively, $A_0=C_0$) too.

\begin{corollary}\label{c:nonsingularity}
Consider the Hermite weight on the real line. With the above notations we have $C_k(\xx)$ nonsingular for all $\xx \in S$.
For $\xx \in X$ the matrices $A_k(\xx)$ are all nonsingular, too.

For $\xx \in X_+\setminus X_-$ only $A_{n+1}(\xx)$ is nonsingular, whereas the remaining $A_k(\xx)$ are singular for $k=0,1,\ldots,n$.

For $\xx \in X_-\setminus X_+$ only $A_{0}(\xx)$ is nonsingular, whereas the remaining $A_k(\xx)$ are singular for $k=1,\ldots,n, n+1$.

For $\xx\in X_-\cap X_+$ the matrices $A_k(\xx)$ are singular for all $k=0,1,\dots,n+1$.
\end{corollary}

\section{The Minimax Equioscillation Theorem}\label{sec:mnimax}

In the cases of the classical and exponential weights, we relied on the fact that the ``difference map'' $\Gamma: S\to \RR^{n+1}$ is a homeomorphism, where
\begin{equation}\label{gamma}
\Gamma(\xx):=(m_1(\xx)-m_0(\xx), \dots, m_{n+1}(\xx)-m_{n}(\xx)).
\end{equation}

However, in the present setting, this property no longer holds. Nevertheless, we can still prove that $\Gamma$ and also $\mm$ defined by
\begin{equation}\label{def_mm}
\mm(\xx):=(m_0(\xx),m_1(\xx),\dots,m_{n+1}(\xx)).
\end{equation}
are proper. Recall that a continuous map $f: X\to Y$ between topological spaces is called \emph{proper}, if $f^{-1}(Q)(\subset X)$ is compact for any compact subset of $Y$.

\begin{proposition}\label{prop:properness}
Let $w$ be a continuous, positive weight function on $\RR$ with $w(t) t^n \to 0$. 
The mappings $\Gamma: S\to \RR^{n+1}$ and $\mm: S\to \RR^{n+2}$ are proper on $S$.
\end{proposition}

This property for the two mappings, with minor differences in the settings, has also been established in earlier works; see \cite{CBoorPinkus, SzAMH}. The proof given in \cite{Exp} can be adapted to the present setting with only minor modifications. Hence, we omit the full proof and only indicate the necessary changes to adapt the argument.

Indeed, the proof of Lemma 15 in \cite{Exp} can be adapted almost verbatim. The only essential difference is that \cite{Exp} works on the positive half-line, while in the present setting one must also treat the possibility that $x_0\to-\infty$. The corresponding argument is symmetric to the case $x_n\to+\infty$ considered in \cite{Exp}$,$ and therefore the properness of $\Gamma$ follows in exactly the same manner.

Following Corollary 7 of \cite{Exp}, we also obtain the properness of the mapping $\mm$.

\medskip

We now turn to the minimax problem of finding optimal interpolation nodes. More precisely, we consider
\begin{align*}
\inf_{\xx \in S} \overline{m}(\xx):=\inf_{\xx \in S} \max_{i=0,\ldots,n+1} m_i(\xx)
= \inf_{\xx \in S} \sup_\RR L_n(\xx,\cdot) = \inf_{\xx \in S} \|L_n(w,\xx,\cdot)\|.
\end{align*}

\begin{theorem}\label{th:minimax} Consider the Hermite weight $w(t):=\exp(-t^2)$, and the weighted Lagrange interpolation in the space $\WW_n$.

There exists an extremal minimax node system $\yy \in S$ such that $\|L_n(w,\yy,\cdot)\| \le \|L_n(w,\xx,\cdot)\|$ for any other node system $\xx \in S$.

Moreover, any such node system belongs to $X$ and equioscillates: $m_0(\yy)=\dots=m_{n+1}(\yy)$.
\end{theorem}

\begin{proof} The first part--existence of a minimax node system in $S$--follows from continuity and properness of the interval maxima mapping $\mm$, see Proposition \ref{prop:properness}.

First, we claim that any minimax node system necessarily belongs to $X$. In other words, we are to prove that a point $\yy \in X^c=S\setminus X$ cannot be a minimax point system. Indeed, if $\yy \in X^c$, then either there is a unique interval maximum equal to 1 ($m_0(\yy)$ or $m_{n+1}(\yy)$), or both of these two interval maxima are equal to 1, while all other interval maxima $m_i(\yy)>1$, $(i=1,\ldots,n)$, as was shown using strict log-concavity, see Lemma \ref{l:strict_concave_place_zi}.

Assume first that $m_{n+1}(\yy)=1$, i.e., $\yy=(y_0,\ldots,y_n)\in X_+$. Since $m_i(\yy)>1$ for $i=1,\ldots,n$, there exists $\ve>0$ such that $m_i(\yy)>1+\ve$, $i=1,\ldots,n$. If, in addition, $\yy\in X_+\setminus X_-$, then $\ve$ may be chosen so that $m_0(\yy)>1+\ve$ as well.

By continuity, there exists $\de>0$ such that, whenever $\|\xx-\yy\|<\de$, $m_i(\xx)>1+\frac{\ve}{2}$,
$i=1,\ldots,n$, and $m_{n+1}(\xx)<1+\frac{\ve}{2}$. If $\yy\notin X_-$, we may also require $m_0(\xx)>1+\frac{\ve}{2}$. Since $\yy\in X_+$, we can apply the derivative properties \eqref{signchange_mi} and \eqref{signchangeXc}, i.e.,
$$\partial m_i(\yy)/\partial x_n>0, \quad i=1,\dots, n; \qquad \partial m_0(\yy)/\partial x_n\ge 0, \qquad \partial m_{n+1}(\yy)/\partial x_n= 0.$$
Consequently, decreasing $y_n$ to some $x_n\in(y_n-\delta,y_n)$ sufficiently close to $y_n$, the resulting node system $\xx=(y_0,\dots, y_{n-1},x_n)$
satisfies 
\begin{align}\label{eq:minimax}
m_0(\xx)&\le m_0(\yy); \nonumber\\
m_i(\xx)&<m_i(\yy), \qquad \quad \text{ for every } i=1,\dots,n;\\
(1=m_{n+1}(\yy)\le) m_{n+1}(\xx)&<1+\ve/2. \nonumber
\end{align}

Now suppose that $\yy\in X_+\setminus X_-$. Then $m_0(\yy)>1$. Let $M$ denote the operator norm corresponding to $\yy$, i.e., $M:=\|L_n(w,\yy,\cdot)\|$. Since $M\ge m_i(\yy)>1+\varepsilon$, $i=1,\ldots,n$, we have
$m_{n+1}(\xx)<1+\frac{\varepsilon}{2}<M$.

If $m_0(\xx)<m_0(\yy)$, then $m_i(\xx)<m_i(\yy)\le M$, $i=0,\ldots,n$, and therefore $\|L_n(w,\xx,\cdot)\|<M$.

Otherwise, suppose that $m_0(\xx)=m_0(\yy)>1+\varepsilon$. If $m_0(\yy)<M$, then all interval maxima corresponding to $\xx$ are still strictly smaller than $M$, and hence $|L_n(w,\xx,\cdot)|<M$.
However, if $m_0(\xx)=m_0(\yy)= \|L_n(w,\yy,\cdot)\|$, then we get only $\|L_n(w,\xx,\cdot)\| =\|L_n(w,\yy,\cdot)\|$. So, we need a further modification of $\xx$. Applying the continuity of the functions $m_i$ once more, we may decrease only the first coordinate of $\xx$, obtaining a node system $\xx'=(x_0',y_1,\ldots,y_{n-1},x_n)$, where $x_0'<y_0$ is chosen sufficiently close to $y_0$, such that $m_i(\xx')< m_i(\yy)$ for $i=1,\ldots, n$. Moreover, $m_{n+1}(\xx')<1+\ve<M$ and the inequality \eqref{signchangeX+} $(\partial m_0(\xx)/\partial x_0>0$ if $\xx\in S \setminus X_-)$ implies, that $m_0(\xx')<m_0(\xx)=m_0(\yy)$. Hence $m_i(\xx')<M$ for every $i=0,\ldots,n+1$ and therefore $\|L_n(w,\xx',\cdot)\| < \|L_n(w,\yy,\cdot)\|$.

Now, let $\yy\in X_+\cap X_-$, i.e., $z_0(\yy)=y_0$, then $m_0(\yy)=1$. Since $m_i(\yy)>1$ for $i=1,\ldots,n$, it implies, that $m_0(\yy)\neq \|L_n(w,\yy,\cdot)\|$. Applying the same properties \eqref{signchange_mi} and \eqref{signchangeXc} for the partial derivatives $\partial m_i/\partial x_n(\yy)$, we get that the second and third inequalities in \eqref{eq:minimax} satisfy, but concerning $m_0$, we get that $m_0(\xx)\le m_0(\yy)=1$ (implying that $m_0(\xx)=1$). Then, we again obtain, that $\|L_n(w,\xx,\cdot)\| < \|L_n(w,\yy,\cdot)\|=M$, and hence $\yy$ cannot be a minimax node system.

Consequently, no node system in $X_+$ can be minimax. By symmetry, no node system in $X_-$ can be minimax either. Hence every minimax node system belongs to $X$.

However, on $X$ we have  $D_k\ne 0$ for all $k$,
in view of Corollary \ref{c:nonsingularity}.
From here, we will follow the direct perturbation argument of Kilgore \cite{Kilgore} and Carl de Boor and Pinkus \cite{CBoorPinkus}.

Assume $\uu\in X$ satisfying $m_k(\uu)< \|L_n(w,\uu,\cdot)\|$ for some index $k$. Let $\ve= \|L_n(w,\uu,\cdot)\|-m_k(\uu)$. By the continuity of the functions $m_i$ we get that there exists $\de>0$ such that for any $\xx \in S$ with $\|\xx-\uu\|<\de$ we have $|m_i(\xx)-m_i(\uu)|<\ve/2$ ($i=0,\ldots,n+1$).
Since $X$ is open, we have
$\uu\in X$.
Moreover, $D_{k}(\uu)\ne 0$, that is $A_{k}(\uu)$ is nonsingular, which entails that in a sufficiently small neighborhood of $\uu$ the mapping $\xx \to (m_0(\xx), \ldots, m_{k-1}(\xx), m_{k+1}(\xx),\ldots,m_{n+1}(\xx))$ acts as a homeomorphism.
So, prescribing the values $\mu_i < m_i(\uu)$ ($i=0,\ldots,n+1$, $i\ne k$) sufficiently close to the original values of $m_i(\uu)$, there is a node system $\vv \in S$,
(actually, $\vv\in X$)
close to $\uu$ (in particular closer than $\de$) with these prescribed values, i.e., with $m_i(\vv)=\mu_i$. As $\| \uu-\vv\|<\de$, the value $m_{k}(\vv)$ also satisfies  $m_{k}(\vv)<m_{k}(\uu)+\ve/2 < \|L_n(w,\uu,\cdot)\|$. It follows that  $\|L_n(w,\vv,\cdot)\| <  \|L_n(w,\uu,\cdot)\|$. Therefore, $\uu$ is not a minimax node system.

Consequently, for a minimax point $\yy$ it is necessary that all values $m_i(\yy)$ coincide, i.e., $\yy$ is an equioscillating node system.
\end{proof}

In fact, more can be said than in the previous Theorem \ref{th:minimax}, as we present below. Namely, for every equioscillating node system $\yy$, the sequence $D_i(\yy)$ has alternating signs. This was already shown by Carl de Boor and Pinkus \cite{CBoorPinkus}, and an elegant proof based on a linear programming theorem was later given by Shi in Lemma 2 of \cite{Shi}. In \cite{Exp}, we further proved that this alternating property remains valid even if the condition $D_i(\xx) \ne 0$ is assumed only on an open subset $Y \subset S$.

\begin{lemma}\label{l:alternate} Let $m_i:S\to \RR$ $(i=1,\ldots,n+1)$ be nonnegative, continuously differentiable functions and let $Y \subset S$ be any open subset of $S$. Assume that $D_k(\xx)\not=0$ for every $\xx\in Y$
and $k=0,1,\ldots,n+1$.

Then, if $\yy\in Y$ is a (locally) minimax node system, then it must be equioscillating.

Moreover, the signs of the nonzero Jacobian determinants  satisfy $(-1)^{k-i}D_i(\yy)D_k(\yy)>0$.
\end{lemma}

\begin{proof}
See Lemma 15 in \cite{Exp}.
\end{proof}

\section{Unicity of the minimax node system}\label{sec:unicity}

We would like to prove the unicity of the found minimax node system. In earlier works \cite{CBoorPinkus, SzAMH, Exp} this always followed from a stronger property, namely the homeomorphism property of the difference mapping $\Gamma$ on $S$. However, in the Hermite-weighted case, $\Gamma$ is no longer a homeomorphism. Indeed, in the matrix  $A$, whose entries are the partial derivatives $\partial m_i/\partial x_j(\xx)$, it may occur in the Hermite case that two column vectors are zero. Therefore, the dimension of the subspace spanned by the nonzero column vectors is at most $n$. Consequently, if we consider any submatrix $A_k$ obtained by selecting an arbitrary $n+1$ columns from $A$, then its determinant is necessarily zero.

This observation shows that the homeomorphism approach cannot be carried out on the whole set S. We therefore restrict our attention to finding a homeomorphism on a suitable subset of X. However, the mapping $\Gamma$ itself is still not suitable for this purpose, since it is not proper on X: as one approaches the boundary between X and its complement $X^c=S\setminus X$, the interval maxima $m_0$ and $m_{n+1}$ tend to $1$. We therefore seek a modified mapping whose components diverge to $\pm\infty$ at the boundary instead of converging to 1. Hence, we introduce a ``modified difference map'' $\Gamma_1$ on $S\setminus X_-$ in the following way.

\begin{definition} Consider the Hermite weight $w(t):=\exp(-t^2)$, and the weighted Lagrange interpolation in the space $\WW_n$ on the real line.

With the above notations, the map $\Gamma_1\colon (S\setminus X_-)\to \RR^{n+1}$ is defined by
\begin{equation}\label{gamma_1}
\Gamma_1(\xx):=(m_1(\xx)-\mu_0(\xx), m_2(\xx)-m_1(\xx), \dots, m_{n+1}(\xx)-m_{n}(\xx)),
\end{equation}
where $\mu_0(\xx):=m_0(\xx)+\log\left(1-\frac{1}{m_0(\xx)}\right)$.
\end{definition}

Our goal is to prove that $\Gamma_1$ is a homeomorphism. To this end we first establish the connectedness of $S\setminus X_-$, and then prove that $\Gamma_1$ is proper.

We begin by proving the connectedness of $X$, for which two technical lemmas are needed. This property was already a key step in the study of exponentially weighted Lagrange interpolation, as shown in \cite{Exp}. However, the structure of $X$ in that setting was somewhat different (in particular, only the condition $m_{n+1}(\xx)=1$ could occur), and therefore a new proof is required in the present framework.

\begin{lemma}\label{l:connected_pre}
Consider the Hermite weight on the real line. With the above notations the following hold.
\begin{itemize}
\item[a)]
For every $a$ there exists $\delta(a)>0$ with the following property. Whenever
\[
a\le u_1<u_2<\dots<u_{n-1}< u_n,
\]
then
\[
\uu_h^{(1)}:=(u_1-h,u_1,\dots,u_{n-1},u_n),
\]
belongs to $S\setminus X_-$ for every $0<h\le \delta(a)$.

\item[b)]
For every $b$ there exists $\delta(b)>0$ with the following property. Whenever
\[
u_0<u_1<u_2<\dots<u_{n-1}< u_n\le b,
\]
then
\[
\uu_h^{(0)}:=(u_0,u_1,\dots,u_{n-1},u_{n-1}+h),
\]
belongs to $S\setminus X_+$ for every $0<h\le \delta(b)$.

\item[c)]
For every $a<b$ there exists $\delta(a,b)>0$ with the following property. Whenever
\[
a\le u_1<u_2<\dots<u_{n-1}\le b,
\]
then
\[
\uu_h:=(u_1-h,u_1,\dots,u_{n-1},u_{n-1}+h)
\]
belongs to $X$ for every $0<h\le \delta(a,b)$.

\end{itemize}
\end{lemma}

\begin{proof}
We first prove Part (a). Consider the weighted polynomials $P_0$ and $P_{n+1}=(-1)^n P_0$ for $\uu_h^{(1)}$. Then, we have $(P_{0}'/P_{0})(t) = \sum_{j=1}^n \frac{1}{t-\eta_j}+\frac{w'(t)}{w(t)}$, where $\eta_j \in I_j$ ($j=1,\dots,n$) are the zeroes of $P_0$. All of them are greater than  $u_1-h$. Taking $t=u_1-h$, all terms of the form $1/(t-\eta_j)$ are negative, so dropping them for $2\le j \le n$ can only increase its value.
Applying that $w$ is log-concave, i.e., $\omega:=\log w$ concave, $0< h \le \delta(a)$ and $P_0(\uu_h^{(1)}, u_1-h)=1$, it follows that \[ P_0'(\uu_h^{(1)}, u_1-h) = \frac{P_0'(\uu_h^{(1)}, u_1-h)}{P_0(\uu_h^{(1)}, u_1-h)} <\frac{w'(u_1-h)}{w(u_1-h)}+\frac{1}{(u_1-h)-\eta_1}\le \omega'(u_1-\delta(a))-\frac{1}{\delta(a)}. \] Since our aim is to find an $\uu_h^{(1)}$ belonging to $S\setminus X_-$, we need to ensure that $P_0'(\uu_h^{(1)}, u_1-h)<0$.
Therefore, if
\begin{equation}\label{delta_assump_1}
\omega'(a-\delta(a))-\frac{1}{\delta(a)}<0,
\end{equation}
then $P_0'(\uu_h^{(1)},u_1-h)<0$ and and the choice of $\delta(a)$ is independent of the node points $u_1,\dots, u_n$.
Consequently, $m_0(\uu_h^{(1)})>1$, which proves part (a).

\medskip

Next, we prove part (b).
This is quite symmetric to the previous case, but we write down its proof for later reference to the derivative estimate occurring below. Consider now the weighted polynomial
$
P_{n+1}=(-1)^n P_0
$
corresponding to
$\uu_h^{(0)}$.
Since $u_{n-1}+h\le u_{n-1}+\delta(b)$, we get that
\[ P_{n+1}'(\uu_h^{(0)},u_{n-1}+h)> \omega(u_{n-1}+\delta(b))+\frac{1}{u_{n-1}+h-u_{n-1}}\ge \omega'(u_{n-1}+\delta(b))+\frac{1}{\delta(b)}.
\]
In order to deduce $\uu_h^{(0)}\in S\setminus X_+$, it suffices to have $P_{n+1}'(\uu_h^{(0)},u_{n-1}+h)>0$.  Consequently, we obtain that if
\begin{equation}\label{delta_assump}
 \omega'(b+\delta(b))+\frac{1}{\delta(b)}>0,
\end{equation}
then $m_{n+1}(\uu_h^{(0)})>1$, i.e., $\uu_h^{(0)}\in S\setminus  X_{+}$
for every $0<h\le\delta(b)$.

Finally, part (c) follows by combining parts (a) and (b). Indeed, choosing
\[
\delta(a,b):=\min\{\delta(a),\delta(b)\},
\]
we obtain simultaneously \eqref{delta_assump_1} and \eqref{delta_assump} as upper resp. lower estimates of the derivatives of $P_0$ resp. $P_{n+1}$, ensuring that the node system is neither in $X_{-}$ nor in $X_{+}$, as before.
\end{proof}

\begin{lemma}\label{l:Xconnected_hermite}
Consider the Hermite weight on the real line.
With the above notations the set $X=\{\xx\in S
z_0(\xx)<x_0<\dots<x_n<z_{n+1}(\xx)\}$ is a domain (i.e., an open connected set) in $\RR^{n+1}$.
\end{lemma}

\begin{proof}
We already know that $X$ is open, for it is an intersection of two strict level sets of continuous functions, namely $X=\{\xx\in S: z_0(\xx)<x_0\}\cap \{\xx\in S: z_{n+1}(\xx)>x_n\}$.

Let $\xx=(x_0,\dots,x_n)$ and $\yy=(y_0,\dots,y_n)$ be two arbitrary elements of $X$. Our aim is to connect them with a curve in $X$. Let $a:=\min\{x_1,y_1\}$ and $b:=\max\{x_{n-1},y_{n-1}\}$ and let us apply the previous Lemma \ref{l:connected_pre}. Pick the appropriate $\delta>0$, furnished by the Lemma. Then we know that points $\uu \in S $ satisfying $a \le u_1 < \dots <u_{n-1}\le b$ have the property that replacing $u_1-h$ for $u_0$ and $u_{n-1}+h$ for $u_n$, the arising $\uu_h$ belong to $X$, whenever $0<h<\delta$. In particular, if $h <\de$, then also $\xx_h, \yy_h \in X$.

So here we choose some $h \in (0,\de)$ such that we also have $h < \min(x_1-x_0, x_n-x_{n-1}, y_1-y_0, y_n-y_{n-1})$. Then we fix this value of $h$ and the respective auxiliary points $\xx_h, \yy_h \in X$ for the rest of the argument. Note that $x_0<x_1-h$, $x_{n-1}+h<x_n$, and similarly, $y_0<y_1-h$, $y_{n-1}+h<y_n$ by construction.

We will then connect $\xx$ to $\yy$ within $X$ through these auxiliary points $\xx_h$ and $\yy_h$. The argument can be divided into three parts. Namely, first we connect $\xx$ to $\xx_h:=(x_1-h,\xx',x_{n-1}+h)$, where $\xx':=(x_1,\dots,x_{n-1})$. Quite similarly, we are to connect $\yy_h$ to $\yy$. Finally, $\xx_h$ will be connected to $\yy_h:=(y_1-h,\yy',y_{n-1}+h)$.

Connecting $\xx$ and $\xx_h$ also contains two steps, the first being the connection of $\xx$ to $\xx_h^*:=(x_1-h,x_1,x_2,\ldots,x_{n-1},x_n)$. We have chosen $h$ to be less than $x_1-x_0$, so that $x_0<x_1-h$ holds. Then the connecting curve is just the segment $[\xx,\xx_h^*]$, so that all coordinate values $x_i$ ($i=1,\ldots,n$) remain fixed, and only the 0th coordinate changes.
Applying the inequalities \eqref{signchangeX+} and \eqref{signchange_outer_m}  ($\partial m_0(\xx)/\partial x_0>0$ and $\partial m_{n+1}(\xx)/\partial x_0>0$) we get that the change of the values of $m_0$ and $m_{n+1}$ along the way from $\xx$ to any point of the segment $[\xx,\xx_h^*]$ is always positive, so that together with $\xx$ also each points of the segment $[\xx,\xx_h^*]$ belong to $X$.

Similarly, leaving $x_0'=x_1-h$, $x_1,\dots, x_{n-1}$ invariant, we move $x_n$ to $x_{n-1}+h$, again taking into account $x_{n-1}+h<x_n$. Now the connecting curve is the segment $[\xx_h^*,\xx_h]$. Using  the sign of partial derivatives \eqref{signchange_outer_m} and \eqref{signchangeX-} ($\partial m_{0}(\xx)/\partial x_n< 0$ and $\partial m_{n+1}(\xx)/\partial x_n<0$) will then furnish that $[\xx_h^*,\xx_h] \subset X$.

In a similar way, $\yy_h$ and $\yy$ can be connected within $X$ by two segments.

Finally we show, that $\xx_h$ and $\yy_h$ can be connected in $X$, more precisely, that the segment $[\xx_h,\yy_h]$ fully belongs to $X$. Let us introduce the parameterized curve $\uu(s):=(1-s)\xx_h+s\yy_h$, $s\in [0,1]$. Observe that by construction $u_1(s)-u_0(s)=u_n(s)-u_{n-1}(s)=h$ for all $0\le s \le 1$, whereas $u_i(s)=(1-s)(x_h)_i+s(y_h)_i \in [a,b]$. So, Lemma \ref{l:connected_pre} applies and $\uu(s) \in X$. Thus we get that $\xx_h$ and $\yy_h$ are connected within $X$ with the curve $\uu(s)=(1-s)\xx_{h}+s\yy_{h}$.
\end{proof}

Moreover, we need to prove, that $S\setminus X_-$ is also connected.

\begin{lemma}\label{l:X_1 domain}
Consider the Hermite weight on the real line. With the above notations the set $S\setminus X_-$ is a domain in $S$ and hence in $\RR^{n+1}$.
\end{lemma}

\begin{proof}
Since $S\setminus X_-$ is a strict level set of continuous function, namely $S\setminus X_-=\{\xx\in S: z_0(\xx)<x_0\}$, it is an open set.

\medskip

Next, we prove connectedness.

Since $X$ is connected by Lemma \ref{l:Xconnected_hermite}, and $S\setminus X_-=X\cup (X_+\setminus X_-)$, it is enough to show that every point of $X_+\setminus X_-$ can be connected to a point of $X$ by a continuous curve contained in $S\setminus X_-$.

Let $\xx=(x_0,\dots,x_n)\in X_+\setminus X_-$. Then $m_0(\xx)>1$, i.e., $z_0(\xx)<x_0$, and $m_{n+1}(\xx)=1$. Fix the coordinates $x_1,\dots,x_n$ and consider the segment
\[
\xx(s):=(x_0+s,x_1,\dots,x_n),
\qquad
0\le s\le x_1-h-x_0,
\]
where $h>0$ will be chosen later.

Since $\xx\in X_+\setminus X_-\subset S\setminus X_-$, by the inequality \eqref{signchangeX+}, we have
\[
\frac{\partial m_0}{\partial x_0}(\xx)>0,
\]
hence $m_0(\xx(s))>1$ along the whole segment. Since $m_{n+1}(\xx)\ge 1$ holds for any $\xx\in S$, it follows that $\xx(s)\in S\setminus X_-$. In particular, with $s=x_1-h-x_0$ the node system $
\xx^*:=(x_1-h,x_1,\dots,x_n)\in S\setminus X_-$.

Now choose $h>0$ according to Lemma \ref{l:connected_pre}~(c), and denote
\[
\xx_h:=(x_1-h,x_1,\dots,x_{n-1},x_{n-1}+h).
\]
By construction, $\xx_h\in X$. It remains to connect $\xx^*:=(x_1-h,x_1,\dots,x_n)$
to $\xx_h$.

For this, consider the segment
\[
\gamma(t):=
(x_1-h,x_1,\dots,x_{n-1},(1-t)x_n+t(x_{n-1}+h)),
\qquad t\in[0,1].
\]

Along this segment, the first $n$ coordinates remain fixed, while the last coordinate decreases.

Since the first coordinate is fixed equal to $x_1-h$, Lemma \ref{l:connected_pre}~(a) implies that every point of the segment satisfies
\[
m_0>1.
\]
Moreover, always
\[
m_{n+1}\ge1.
\]
Therefore,
\[
\gamma(t)\in S\setminus X_-
\qquad (0\le t\le1).
\]

Since $\gamma(1)=\xx_h\in X$, we have connected $\xx$ to a point of $X$ by a continuous curve contained in $S\setminus X_-$. Thus $S\setminus X_-$ is connected.

Consequently, we get that $S\setminus X_-$ is a domain in $S$.
\end{proof}


\begin{lemma}\label{l:gamma1_proper}

Let $w$ be a continuous, positive weight function on $\RR$ with $w(t)|t|^n\to 0$ ($|t|\to\infty$). 
Then the mapping $\Gamma_1\colon (S\setminus X_-)\to \RR^{n+1}$ is proper.
\end{lemma}

Throughout the remainder of the paper, $\partial_A B$ denotes the boundary of a subset $B\subset A$ relative to the space $A$.

\begin{proof}
Let $Q\subset \RR^{n+1}$ be compact, and define $W:=\Gamma_1^{-1}(Q)$. Since $\Gamma_1$ is continuous, $W$ is relatively closed in $S\setminus X_-$.

Since $Q$ is bounded, there exists a constant $C>0$ such that
\[
|m_1(\xx)-\mu_0(\xx)|\le C,
\qquad
|m_{i+1}(\xx)-m_i(\xx)|\le C,
\quad i=1,\dots,n,
\]
for every $\xx\in W$, where
\[
\mu_0(\xx):=\mu(m_0(\xx)),
\qquad
\mu(t):=t+\log\left(1-\frac{1}{t}\right), \qquad t>1.
\]

Since $m_1(\xx)\ge 1$ and $|m_1(\xx)-\mu_0(\xx)|\le C$, we have $\mu_0(\xx)>1-C$. Since $\mu$ is strictly increasing and $\mu(t)\to -\infty$ as $t\to 1+$, the set $\{t>1:\mu(t)\ge 1-C\}$ is of the form $[1+\ve,\infty)$ for some $\ve>0$. Hence $m_0(\xx)\ge 1+\ve$ for every $\xx \in W$.

Consequently, $0>\mu_0(\xx)-m_0(\xx)=\log\left(1-\frac{1}{m_0(\xx)}\right)\ge \log\left(1-\frac{1}{1+\ve}\right)$ remains bounded on $W$. Hence
$m_1(\xx)-m_0(\xx)=(m_1(\xx)-\mu_0(\xx))+(\mu_0(\xx)-m_0(\xx))$ is bounded on $W$ as well.

Therefore all coordinates of $\Gamma=\bigl(m_1-m_0,m_2-m_1,\ldots,m_{n+1}-m_n)$ are bounded on $W$, and thus $\Gamma(W)$ is bounded. Let $K:=\overline{\Gamma(W)}$. Then $K$ is compact. Since $W\subset \Gamma^{-1}(K)$, and $\Gamma$ is proper by Proposition~\ref{prop:properness}, it follows that $\Gamma^{-1}(K)$ is compact. Consequently $W$ is bounded.

Hence, $W$ is bounded and it is relatively closed in $S\setminus X_-$ by the continuity of $\Gamma_1$. In particular, there is a constant $c$ such that $\|\xx\|\le c$ for all points $\xx \in W$. It remains to see that $W$ is closed in $\RR^n$, too. In other words, we have to show that if a point $\xx \in \partial_{\overline{S}} (S\setminus X_-)$ satisfies $\|\xx\|\le c$, then it cannot be a limit point of $W$.

Now, let $\xx \in \partial_{\overline{S}} (S\setminus X_-)$. Then either $\xx\in \partial_{\overline S}S$, or $\xx\in \partial_S X_-$. In the former case, at least one degeneracy occurs, that is, there exists an index $i\in \{1,\ldots,n\}$ such that $x_{i-1}=x_i$. In the latter case, if no such degeneracy occurs, i.e., $x_{i-1}\neq x_i$ $(i=1,\ldots,n)$, then the boundary is attained through the defining condition of $X_-$, namely $z_0(\xx)=x_0$.

First, assume that $x_{i-1}=x_i$ for some $i\in \{1,\dots,n\}$, and let $|\xx|\le c$. Let $\de>0$ be small, and let $\yy\in X$ satisfy $\|\xx-\yy\|<\de$. Take $z_i\in \intt I_i(\yy)$ such that $m_i(\yy)=P_i(\yy,z_i)$. Choose a point $s\in[0,1]$ whose distance from every node $y_0,\dots,y_n$ is at least $1/(2n+1)$, and let $k$ be the unique index satisfying $s\in I_k(\yy)$. With the unique index $k$ with $s \in I_k$ we obviously have
\[
m_k(\yy) \ge P_k(\yy,s) = \sum_{j=0}^n |h_j(\yy,s)|.
\]
This leads to
\begin{align*}
\frac{|h_j(\yy,s)|}{|h_j(\yy,z_i)|} & = \frac{w(s)}{w(z_i)} \prod_{l=0, l \ne j}^n \left| \frac{s-y_l}{z_i-y_l}\right|
\ge \frac{\min_{[0,1]} w}{\max_{[-c-\de,c+\de]} w} \left( \frac{1}{2n+1} \right)^n (2(c+\de))^{-n+1} \frac{1}{2 \de},
\end{align*}
for the point $z_i$ is in the interval $I_i(\yy)=[y_{i-1},y_i]$ of length at most $2\de$, and at least one of the indices $l=i-1$ or $l=i$ will occur in the product for $l\ne j$, leading to $|z_i-y_l|>1/(2\de)$. This ratio will stay above $C_1/\de$ with some (possibly large) constant $C_1$. Therefore we find
\[
m_k(\yy) \ge \sum_{j=0}^n |h_j(\yy,s)| \ge \frac{C_1}{\de} \sum_{j=0}^n |h_j(\yy,z_i)| = \frac{C_1}{\de} m_i(\yy).
\]
Since $i\in \{1,\ldots,n\}$, we have $m_i(\yy)\ge 1$, and therefore
\[
\max_{l=1,\dots,n+1}|m_l(\yy)-m_{l-1}(\yy)|\ge \frac{1}{n}\left(\frac{C_1}{\de}-1\right).
\]

If this maximum is attained for some index $l\ge 2$, then one of the coordinates of $\Gamma_1(\yy)$ tends to infinity as $\de\to 0$, and thus $\Gamma_1(\yy)\notin Q$ for sufficiently small $\de$.

Assume now that $\max\limits_{l=1,\dots,n+1}|m_l(\yy)-m_{l-1}(\yy)|=|m_1(\yy)-m_0(\yy)|$. Then, $|m_1(\yy)-m_0(\yy)| \to\infty$ if $\de\to 0$. Since $-\log\left(1-\frac{1}{m_0(\yy)}\right)$ is positive and
$$m_1(\yy)-\mu_0(\yy)=m_1(\yy)-m_0(\yy)-\log\left(1-\frac{1}{m_0(\yy)}\right)\ge m_1(\yy)-m_0(\yy),$$
we also have $|m_1(\yy)-\mu_0(\yy)| \to\infty$.

If $m_1(\yy)-m_0(\yy)\to -\infty$, then necessarily $m_0(\yy)\to\infty$, since $m_1(\yy)\ge 1$. Consequently, $\log\left(1-\frac{1}{m_0(\yy)}\right)\to 0$, and therefore
\[
m_1(\yy)-\mu_0(\yy)=m_1(\yy)-m_0(\yy)-\log\left(1-\frac{1}{m_0(\yy)}\right)\to-\infty.
\]

Thus, in all cases, at least one coordinate of $\Gamma_1(\yy)$ diverges, and hence $\Gamma_1(\yy)\notin Q$ for sufficiently small $\de$.

It remains to consider the case when $\xx\in X_-$, while $x_{i-1}\ne x_i$ for all $i=1,\dots,n$.

Applying Lemma~\ref{l:strict_concave_place_zi}, we obtain points $z_i(\xx)\in\intt I_i(\xx)$, for every $i=1,\dots,n$, such that $m_i(\xx)=P_i(\xx,z_i(\xx))>1$, (and $m_{n+1}(\xx)\ge 1$, which is out of interest).

Now, we consider the index $i=1$ and let $$\eta:=m_1(\xx)-1>0.$$ By continuity, there exists $\de>0$ such that $\|\xx-\yy\|<\de$ implies
\[
m_1(\yy)>1+\frac{\eta}{2}.
\]

On the other hand, $m_0(\xx)=1$, and since $m_0$ depends continuously on $\xx$, we have $m_0(\yy)\to 1^+$. Therefore $\mu_0(\yy)=m_0(\yy)+\log\left(1-\frac{1}{m_0(\yy)}\right)\to-\infty$.

Since $m_1(\yy)\to m_1(\xx)>1$, it follows that $m_1(\yy)-\mu_0(\yy)\longrightarrow+\infty$.
Hence, $\|\Gamma_1(\yy)\|\to\infty$, and therefore $\xx$ cannot be a limit point of $W$.

Thus $W$ is closed in $\RR^{n+1}$. Being bounded and closed, it is compact. Consequently $\Gamma_1$ is proper.
\end{proof}

\begin{corollary}\label{c:proper_M1}
Let $w$ be a continuous, positive weight function on $\RR$ with $w(t)|t|^k\to 0$ ($|t|\to\infty$) for all $k$ (or at least for $k=n$). Then, the mapping ${\mathbf M_1}: (S\setminus X_-)\to \RR^{n+2}$ defined by
\begin{equation}\label{M1def}
\mathbf{M_1}(\xx)=(\mu_0(\xx), m_1(\xx),\dots,m_{n+1}(\xx)),
\end{equation}
is a proper mapping from
$\setminus X_{-}$ to $\RR^{n+2}$.
\end{corollary}
\begin{proof}
The proper map $\Gamma_1$ is a composition of the continuous mappings $\mathbf{M_1}: (S\setminus X_-)\to \RR^{n+1}$ and the simple ``difference mapping'' $d: \RR^{n+2}\to \RR^{n+1}$ given by
\[
d(\mu_0, m_1, \dots,m_{n+1}) := (m_1-\mu_0,\dots,m_{n+1}-m_n).
\]
That is, $\Gamma_1(\xx)=(d\circ \mathbf{M_1})(\xx)$. The proof can now be completed by following the argument used in the proof of Corollary 7 in \cite{Exp}.
\end{proof}

\begin{theorem}\label{t:homeom}
Consider the Hermite weight on the real line.
Then, the mapping $\Gamma_1\colon (S\setminus X_-) \to \mathbb{R}^{n+1}$, defined by \eqref{gamma_1}, is a homeomorphism.
\end{theorem}

\begin{proof}
According to Lemma \ref{l:gamma1_proper}, $\Gamma_1$ is a proper mapping such that its image space is $\RR^{n+1}$.
And we know, that $S\setminus X_-$ is a connected open set (an open domain) (see Lemma \ref{l:X_1 domain}).

Then, knowing the properness a well-known theorem of topology going back to Hadamard says that $\Gamma_1$ will be a global homeomorphism, if and only if it is to be a local homeomorphism.
In turn, local homeomorphism would follow, in view of continuous differentiability, from
the nonvanishing of the Jacobian of the derivative matrix
\begin{equation}
D\Gamma_1=\begin{aligned}
\left(\begin{array}{cccc}
\frac{\partial m_{1}}{\partial x_0}(\xx)-\frac{\mu_0(\xx)}{\partial x_0}(\xx)& \frac{\partial m_{1}}{\partial x_1}(\xx)-\frac{\partial \mu_0}{\partial x_1}(\xx)& \ldots & \frac{\partial m_{1}}{\partial x_n} (\xx)-\frac{\partial \mu_0}{\partial x_n}(\xx)\\
\frac{\partial m_{2}}{\partial x_0}(\xx)-\frac{\partial m_1}{\partial x_0}(\xx) & \frac{\partial m_{2}}{\partial x_1}(\xx)-\frac{\partial m_1)}{\partial x_1}(\xx) & \ldots & \frac{\partial m_{2}}{\partial x_n}(\xx)-\frac{\partial m_1}{\partial x_n}(\xx)\\
\vdots & \vdots &  \ddots & \vdots\\
\frac{\partial m_{n+1}}{\partial x_0}(\xx)-\frac{\partial m_n}{\partial x_0}(\xx) & \frac{\partial m_{n+1}}{\partial x_1}(\xx)-\frac{\partial m_n}{\partial x_1}(\xx) & \ldots & \frac{\partial m_{n+1}}{\partial x_n}(\xx)-\frac{\partial m_n}{\partial x_n}(\xx)
\end{array}
\right) .
\end{aligned}
\end{equation}
Firstly, we calculate
\[
\frac{\partial m_1}{\partial x_j}(\xx)-\frac{\partial \mu_0}{\partial x_j}(\xx)=\frac{\partial m_{1}}{\partial x_j}(\xx)-\left(1+\frac{1}{m_0(\xx)(m_{0}(\xx)-1)}\right)\frac{\partial m_{0}}{\partial x_j}(\xx)=:\frac{\partial m_{1}}{\partial x_j}(\xx)-\lambda(\xx)\frac{\partial m_{0}}{\partial x_j}(\xx)
\]
where $$\lambda(\xx):=1+\frac{1}{m_0(\xx)(m_0(\xx)-1)}$$ is positive for all $\xx\in S\setminus X_-$.

Let $\ba_i$ denote the $i$th column of $A$. Then, we can write
$$
\det(D\Gamma_1)=\det[\ba_1-\lambda(\xx)\ba_0, \cdots, \ba_{i+1}- \ba_i, \cdots \ba_{n+1}- \ba_n]
$$
Applying Lemma 16 in \cite{Exp}, using the vectors $\nu_i a_i$ instead of the vectors $a_i$, we obtain that
\begin{equation}\label{DGsignedsum}
\det D\Gamma_1(\xx) = -D_0(\xx)+\lambda(\xx)\sum_{k=1}^{n+1} (-1)^{k+1} D_k(\xx).
\end{equation}

Next we are to show $\det D\Gamma_1 \ne 0$ separately both for $X$ and $X_+\setminus X_-$.
As for $X_+\setminus X_-$, the situation is relatively simple. According to Corollary \ref{c:nonsingularity}, all $A_k$ are singular (with the last column being identically zero) except for $A_{n+1}$, which is not. So for the respective Jacobian determinants, $D_k(\xx)=0$, $k=0,\ldots,n$, while $D_{n+1}(\xx)\ne 0$, whenever $\xx\in X_+\setminus X_-$.
Therefore we have that $\det D\Gamma(\xx)=(-1)^{n+2}\lambda(\xx)D_{n+1}(\xx)\ne 0$ on $X_+\setminus X_-$.

Now, we turn to the case $\xx\in X$.

The remaining arguments follow the classical approach of de Boor and Pinkus \cite{CBoorPinkus}. Restricting ourselves to the regular set $X$, the proof proceeds exactly as in Lemma~16 of \cite{Exp}. In particular, Lemma~\ref{l:alternate} immediately yields the following statement.

\begin{lemma}\label{l:JiJkonX} Let $w(t)=\exp(-t^2)$ and consider the respective Hermite interpolation. With the above notations the following holds. If $\xx\in X$ and $0 \le i,k \le n+1$ then we have
\begin{equation}\label{Di_alternate}
(-1)^{i-k}D_i(\xx)D_k(\xx)>0.
\end{equation}
\end{lemma}

\begin{proof}
According to Lemma \ref{l:alternate}, \eqref{Di_alternate} holds for any minimax node system $\yy\in X$, whereas at least one such minimax node system must exist in $X$ in view of Theorem \ref{th:minimax}. Furthermore, $X$ is connected, according to Lemma \ref{l:Xconnected_hermite} and $D_i(\xx)D_k(\xx)\ne 0$ in $X$ in view of Corollary \ref{c:nonsingularity}. Therefore the expression \eqref{Di_alternate} has the same sign all over $X$.
\end{proof}

Now, obviously equations \eqref{DGsignedsum} and \eqref{Di_alternate} imply that $\det D\Gamma_1 \ne 0$ also holds.

\end{proof}

\begin{corollary}
Let $w(t)=\exp(-t^2)$ and consider the respective Hermite interpolation. Then, the equioscillating node system, and therefore the extremal node system is necessarily unique and belongs to $X$.
\end{corollary}

\begin{proof}
By Theorem \ref{th:minimax}, there exists at least one minimax node system and every minimax node system belongs to $X$, and equioscillates. Therefore, its image under $\Gamma_1$ is the zero vector $\bf{0}$. However, $\Gamma_1$ is a homeomorphism on $X\cup X_1$, hence there exists a unique $\yy\in S\setminus X_-$ such that $\Gamma_1(\yy)=\bf{0}$. Since $m_{n+1}$ is 1 on $X_+\setminus X_-$, hence $\yy\in X$.

Furthermore, if $\xx\in X^c$, then $m_i(\xx)>1$ for every $i=1,\dots,n$, see Corollary \ref{c:ziininttIi} and at least one of the following equations holds: $m_0(\xx)=1$, or $m_{n+1}(\xx)=1$. Consequently, there is no equioscillating node system in $X^c$. Recall that there is neither any extremal node system in $X^c$, according to Theorem \ref{th:minimax}.

Therefore, the unique equioscillating system is the only optimal node system and it belongs to $X$.
\end{proof}

\section{The Erdős Conjecture and issues of majorization}\label{sec:majorization}

\begin{lemma}\label{lem:Cramer} Let $w(t)=\exp(-t^2)$ and consider the respective Hermite interpolation. With the above notations the following holds.

Let $0\le k \le n+1$ and let $\Phi(\xx):=\Phi_k(\xx):=(m_i(\xx))_{i=0,i\ne k}^{n+1}$ ($\xx\in S^{(n+1)}$). Assume that $D_k(\uu)\ne 0$ for a node system $\uu \in S^{(n+1)}$. Then in a sufficiently small neighborhood $U \subset S^{(n+1)}$ of $\uu$, $\Phi:U\leftrightarrow W$ is a homeomorphism between $U$ and $W:=\Phi(U)$, and on $W$ $m_k$ can be expressed as a continuously differentiable function $G=G_k$ of the interval maxima $\ww:=(m_0,\ldots,m_{k-1},m_{k+1},\ldots, m_{n+1})\in W$. Moreover, for $\ww=\Phi(\vv)$ the gradient of $G$ is
\begin{equation}\label{Mnderivative}
\nabla_{\ww} G =
\begin{bmatrix}
. \\ . \\ . \\
\frac{\partial m_k}{\partial m_i}(\ww)
\\ . \\ . \\ .
\end{bmatrix}_{i=0,i\ne k}^{n+1}
= \begin{bmatrix}
. \\ . \\ .\\
\frac{(-1)^{k-i-1}D_i(\vv)}{D_k(\vv)}
\\ . \\ . \\ .
\end{bmatrix}_{i=0, i\ne k}^{n-1}
(\vv \in U, \ww=\Phi(\uu)\in W).
\end{equation}
\end{lemma}
\begin{proof} We already know that all the $m_i$ are continuously differentiable functions of the variables $\xx \in S^{(n+1)}$, hence the function $$\Phi(\xx):= (m_0(\xx),\ldots,m_{k-1}(\xx),m_{k+1}(\xx),\ldots,m_{n+1}(\xx))$$ is continuously differentiable, too. Moreover, wherever $D_k(\xx)\ne 0$, it is a local homeomorphism with derivative matrix $A_k(\xx)$ at $\xx$. Therefore, it has a continuously differentiable inverse function $\Phi^{-1}$, too. So, we can write $m_k(\xx)=m_k(\Phi^{-1}(m_0(\xx),\ldots,m_{k-1}(\xx),m_{k+1}(\xx),\ldots,m_{n+1}(\xx)))$, and $G=m_k \circ \Phi^{-1}$ or $G \circ \Phi= m_k$. From this it is clear that $G$ is a continuously differentiable function of the values $(m_0,\ldots,m_{k-1},m_{k+1},\ldots,m_{n+1})$ in a small neighborhood of any given point $\Phi(\uu)=(m_i(\uu))_{i=0, i\ne k}^{n+1}$ with $D_k(\uu)\ne 0$. Equivalently, for a small neighborhood $U$ of $\uu \in S^{(n+1)}$, $G$ is a continuously differentiable mapping in $W=\Phi(U)$.

So let us differentiate at some $\vv \in U$ the composite function $m_k=G\circ \Phi$. By an application of the chain rule we get for each $j=0,\ldots,n$ that
$$
\frac{\partial m_k}{\partial x_j}(\vv) = \sum_{i=0, i\ne k}^{n+1} \frac{\partial m_k}{\partial m_i}(\Phi(\vv)) \frac{\partial m_i}{\partial x_j}(\vv),
$$
or, in matrix form
\begin{equation}\label{Mnderivative_2}
\begin{bmatrix}
. \\ . \\ . \\
\frac{\partial m_k}{\partial x_j}(\vv)
\\ . \\ . \\ .
\end{bmatrix}_{j=0}^n
= \begin{bmatrix}
. & .& . \\ . & .& .  \\ . & .& . \\
. &  \frac{\partial m_i}{\partial x_j}(\vv) & . \\
. & .& . \\ . & .& .  \\ . & .& .
\end{bmatrix}_{j=0,i=0, i\ne k}^{n,n+1}
\cdot \quad
\begin{bmatrix}
. \\ . \\ . \\
\frac{\partial m_k}{\partial m_i}(\Phi(\vv))
\\ . \\ . \\ .
\end{bmatrix}_{i=0, i\ne k}^{n+1}
.
\end{equation}
The square matrix on the right hand side is $A_k(\vv)$, with determinant $D_k(\vv)\ne 0$. So, this system of equations is a linear equation in $\RR^n$ with nonsingular coefficient matrix $A_k$ and constants $\beta_{j}:=\frac{\partial m_k}{\partial x_j}(\vv)$, considering the partial derivatives $\xi_i:=\frac{\partial m_k}{\partial m_i}(\ww)=\frac{\partial m_k}{\partial m_i}(\Phi(\vv))$ as the $n$ unknowns.
By Cramer's rule, the system can be solved and the unknowns can be represented as
$$
\xi_i=\frac{\det B_i}{D_k},
$$
where $B_i$ is the matrix obtained from $A_k$ by replacing its $i$th column by the vector of constants $\bf{\beta}$. However, this new $i$th column is just the $k$th column of $A$, and exchanging it one by one with the $i+1$st, $i+2$nd etc. $k-1$th columns if $k>i$ and $i-1$st, $i-2$nd etc. $k+1$th if $k<i$, we are led to $A_i$, the matrix obtained from $A$ by deleting its $i$th column. So, $\det B_i=(-1)^{k-i-1} \det A_i=(-1)^{k-i-1} D_i$.
\end{proof}

\begin{corollary}\label{l:pert_mkmi}
Let $w(t)=\exp(-t^2)$ and consider the respective Hermite interpolation. With the above notations the following holds. Let $\xx\in X$. Then, for any fixed $k$ with $0\le k \le n+1$, in a small neighborhood $W$ of $(m_0(\xx),\ldots,m_{k-1}(\xx),m_{k+1}(\xx),\ldots,m_{n+1}(\xx))$ the interval maxima $m_k(\xx)$ can be expressed as a differentiable function of $\ww=(m_1,\ldots,m_{k-1},m_{k+1},\ldots,m_{n+1})$, and $\partial m_k/ \partial m_i<0$ for every $i=0,\dots,n+1$, $i\ne k$.
\end{corollary}

\begin{proof}
We have just seen that $\frac{\partial m_k}{\partial m_i}(\ww) = (-1)^{k-i-1} D_i(\vv)/D_k(\vv)$, if $\vv$ belongs to a small neighborhood of $\xx$ and $\ww=\Phi(\vv)$. Comparing to Lemma \ref{l:JiJkonX} immediately furnishes the negativity of this expression.
\end{proof}

\begin{theorem}[Limitations to majorization]\label{th:nonmajorization}
Let $w(t)=\exp(-t^2)$ and consider the respective Hermite interpolation. With the above notations the following holds.
If for any pair $\xx, \yy \in S$ we have that $m_i(\xx) \le m_i(\yy)$, then either $\xx=\yy$, or both $\xx, \yy \in X^c$.

Such majorizing pairs in $X^c$ do occur, however, even in the following strict sense: it may well happen that $m_i(\xx)<m_i(\yy)$ for $i=0, 1, \ldots,n$, whereas $m_{n+1}(\xx)=m_{n+1}(\yy)=1$ or similarly, $m_i(\xx)<m_i(\yy)$ for $i=1, \ldots,n, n+1$, whereas $m_{0}(\xx)=m_0(\yy)=1$.
\end{theorem}

\begin{remark}
Theorem \ref{th:nonmajorization} immediately yields the following refined monotonicity property. For every $\xx\in X_+$, there exists a node system $\xx'$ such that
\[
m_0(\xx)\le m_0(\xx'), \qquad
m_i(\xx)<m_i(\xx') \quad (i=1,\ldots,n),
\]
and
\[
m_{n+1}(\xx)=m_{n+1}(\xx')=1.
\]
Thus, $\xx'$ strictly majorizes $\xx$ on all intervals between consecutive nodes, while equality is preserved on the last interval $[x_n,\infty)$.

Similarly, for every $\xx\in X_-$, there exists a node system $\xx'$ such that
\[
m_0(\xx)=m_0(\xx')=1,\qquad
m_i(\xx)<m_i(\xx') \quad (i=1,\ldots,n), \quad \text{ and } \quad
m_{n+1}(\xx)\le m_{n+1}(\xx').
\]
Hence, $\xx'$ strictly majorizes $\xx$ on all interior intervals, while equality is preserved on the first interval $(-\infty,x_0]$.

Moreover, these statements can be sharpened in the interior of the corresponding regions. If $\xx\in \operatorname{int}(X_+)\setminus X_-$, then $m_i(\xx)<m_i(\xx') \qquad (i=0,\ldots,n)$,
while $m_{n+1}(\xx)=m_{n+1}(\xx')=1$. Likewise, if $\xx\in \operatorname{int}(X_-)\setminus X_+$,
then $m_i(\xx)<m_i(\xx')$ for $i=1,\ldots,n+1$, while $m_0(\xx)=m_0(\xx')=1$.

In other words, away from the boundary between $X_+$ and $X_-$, the theorem yields a node system that strictly increases all interval maxima except for the distinguished endpoint interval, where equality is retained.

\end{remark}

\begin{proof}[Proof of Theorem \ref{th:nonmajorization}]
Let $\xx=(x_0,\dots,x_n)$ belong to $\xx\in \intt(X_+)\setminus X_-$, i.e., $z_{n+1}(\xx)=x_n$ and $z_0(\xx)<x_0$. Then with any $x_n'>x_n$, the node system  $\xx'=(x_0,\dots,x_{n-1},x_n')$ belongs to $ X_+$, too. In particular $m_{n+1}(\xx')=m_{n+1}(\xx)=1$. Applying the properties of partial derivatives \eqref{signchange_mi} and Remark~\ref{r:strict_monotonicity}, we know that
\begin{align*}
\frac{\partial m_i}{\partial x_n}(\xx)>0 \quad (i=0,1,\dots,n).
\end{align*}
The signs of the above partial derivatives entail that $m_i(\xx') > m_i(\xx)$ for $i=0,\dots,n$, whereas $m_{n+1}(\xx)=m_{n+1}(\xx')=1$. Note that we have constructed a majorant for every $\xx \in \intt(X_+)\setminus X_-$ with the same first $n$ coordinates $x_0,\ldots,x_{n-1}$.

Similarly, if $\xx\in \intt(X_-)\setminus X_+$, then decreasing the first coordinate $x_0$ of $\xx$, we can get a node system $\xx'$, such that $m_{0}(\xx')=m_0(\xx)=1$, $m_i(\xx'')>m_i(\xx)$ $i=1,\dots,n+1$.

The second part of the assertion follows.

Now, assume that $m_i(\xx) \le m_i(\yy)$ such that $\xx$ and $\yy$ do not belong to $X^c=X_-\cup X_+$, simultaneously. It implies that $\yy\in X$, (otherwise $1=m_{n+1}(\yy)\ge m_{n+1}(\xx)\ge 1$ or $1=m_{0}(\yy)\ge m_{0}(\xx)\ge 1$ shows that $\xx\in X^c$, too).

Let $\yy\in X$.
Next we prove that then either $\xx=\yy$, and we have equalities for every $i$, or there exists some strictly majorizing $\yy'\in X$ of $\xx$. First, if $m_i(\yy)=m_i(\xx)$ for all indices $i$, then $m_0(\xx)>1$ and $m_{n+1}(\xx)>1$, hence $\xx\in X\subset S\setminus X_-$. Moreover, $\Gamma_1(\xx)=\Gamma_1(\yy)$, hence by Theorem \ref{t:homeom} we must have $\xx=\yy$, too. Second, if there exists an index $k$ with $m_k(\yy)>m_k(\xx)$, then referring to $D_k(\yy)\ne 0$ (for $\yy \in X$), we can slightly perturb the node system $\yy$ to some $\yy'$ with all values $m_i(\yy')$ with $i\ne k$ prescribed to strictly exceed $m_i(\xx)$ (hence also 1), whereas  $m_k(\yy')$ can still be kept greater than $m_k(\xx)$ (hence also 1).

Moreover, we can assume, that $\xx:=(x_0,x_1,\dots,x_{n-1},x_n)\in S\setminus X_-$. Indeed, if $\xx \in X_-$, then $m_0(\xx)=1$. However, by Part a) of Lemma \ref{l:connected_pre}, there exists $x_0'$ with $x_0<x_0'<x_1$ such that the corresponding node system $\xx':=(x_0',x_1,\dots, x_{n})$ lies in $\partial_S (S\setminus X_-)$. Since $X_-$ is a relatively closed set, $\partial_S (S\setminus X_-)\subset X_-$,  therefore $\xx'$ satisfies $m_0(\xx')=1$. Recall, that $x_0'=\sup\{y: \forall y'<y \text{ the node system } (y',x_1,\dots,x_n)\in X_-\}$. Applying the properties \eqref{signchange_mi_X-}, and \eqref{signchangeXc_2} of partial derivatives, we know that for $\xx\in X_-$
\begin{align*}
\frac{\partial m_i}{\partial x_0}(\xx)<0 \quad (i=1,\dots, n);\qquad
\frac{\partial m_{n+1}}{\partial x_0}(\xx)\le 0;\qquad
\frac{\partial m_0}{\partial x_0}(\xx)= 0.
\end{align*}

Thus, as the first node is increased from $x_0$ to $x_0'$, none of the quantities $m_i$ increases. Since, by assumption, $m_i(\xx)<m_i(\yy),$ for every $i=0,\ldots,n+1$, it follows that
\[
m_i(\xx')\le m_i(\xx)<m_i(\yy), \qquad i=0,\ldots,n+1.
\]

Now choose $\ve>0$ sufficiently small and put $\xx''=(x_0'+\ve,x_1,\ldots,x_n)$.
Then $\xx''\in S\setminus X_-$. By the continuity of the functions $m_i$, the values $m_i(\xx'')$ can be made arbitrarily close to $m_i(\xx')$. Since all the inequalities $m_i(\xx')<m_i(\yy)$ for every $i=0,\ldots,n+1$ are strict, $\ve$ may be chosen so that
\[
m_i(\xx'')<m_i(\yy), \qquad i=0,\ldots,n+1.
\]
Thus $\xx''$ has the required properties.

So we can as well assume that $\yy\in X$ strictly majorises $\xx\in S\setminus X_-$: $m_i(\yy)>m_i(\xx)$, $i=0,1,\ldots,n+1$. Denote the $i$th component of $\mathbf{M_1}$, already defined in \eqref{M1def}, by $\lambda_i(\xx)$, i.e.,
$$\lambda_0(\xx):=\mu_0(\xx)=m_0(\xx)+\log\left(1-\frac{1}{m_0(\xx)}\right),\quad \text{ and} \quad \lambda_i(\xx):=m_i(\xx) \quad \text{for} \quad i=1,\dots,n+1.$$
Then, denote $B:=B(\yy):= \max\limits_{0\le i\le n+1}(\lambda_i(\yy)-\lambda_i(\xx))$ and $b:=b(\yy):=\min\limits_{0\le i\le n+1}(\lambda_i(\yy)-\lambda_i(\xx))$. Strict majorization together with the monotonicity of the function $\mu(t)=t+\log(1-1/t)$, built in the definition of the modified coordinate $\mu_0$ implies that $b>0$. Also, in general for arbitrary $\ww \in S$ denote
\begin{align*}
b(\ww):=\min\limits_{0\le i\le n+1}(\lambda_i(\ww)-\lambda_i(\xx)); \qquad B(\ww):=\max\limits_{0\le i\le n+1}(\lambda_i(\ww)-\lambda_i(\xx)).
\end{align*}
Let us consider the set of points $W:=\{\ww \in S~:~ b(\ww)\ge b, B(\ww)\le B\}$. As $\yy \in W$, this set is nonempty; by continuity of the functions $m_i$ and $\mu_0$, the intersection of these closed level sets are also closed. Further, if $\ww \in W$, then $\max\{\mu_0(\ww),m_1(\ww),\dots,m_{n+1}(\ww)\}\le \max_i\lambda_i(\xx)+B:=\Lambda(\xx)+B$, hence $W \subset \mathbf{M_1}^{-1}([1,\Lambda(\xx)+B]^{n+1})$. The latter set is a compact subset of $S$, for $\mathbf{M_1}$ is a proper mapping according to Corollary \ref{prop:properness}. Therefore, the closed subset $W$ is also a compact subset of $S$. Also, $W \subset X$ because $b(\ww)\ge b>0$ entails $\mu_0(\ww)>\mu_0(\xx)>-\infty$ (hence $m_0(\ww)>1$) and $m_{n+1}(\ww) >1$. By compactness, the continuous function $b(\ww)$ has a finite maximum value $b^*$ on $W$, and $b$ must attain this maximum on $W$. So let us restrict attention to the maximal subset $Z\subset W$, where $b(\ww)=b^*$. As another closed level set, it is a closed, hence compact subset of $W$, which is not empty (in view of the existence of a $b$-maximal point in $W$). Let us now \emph{minimize} the value of $B(\ww)$ on $Z$. As before, on the compact set $Z$ there must exist some point(s) $\zz \in Z \subset W \subset X$ with $B(\zz)=B^*:=\min\{ B(\ww)~:~ \ww \in Z\}$.

Now we claim that $B(\zz)=b(\zz)$ for this node system. Before proving, let us see that it concludes the proof of Theorem \ref{th:nonmajorization}. Indeed, then we would arrive at a node system $\zz$, such that $\lambda_i(\zz)=\lambda_i(\xx)+b^*$ for all $i=0,1,\ldots,n, n+1$, hence $\Gamma_1(\zz)=\Gamma_1(\xx)$ with $\zz\in X$ and $\xx\in S\setminus X_-$, so that taking into account that $\Gamma_1$ is a homeomorphism, also $\zz=\xx$, which is not possible as we had $b^* \ge b(\yy)>0$.

So it remains to prove $B^*=b^*$, i.e., $B(\zz)=b(\zz)$. Let us assume for a contradiction that $B(\zz)>b(\zz)$. Then write for the index sets of maximal and minimal differences $\I$ and $\JJ$:
\begin{equation}
\begin{aligned}
\I:=\{K~:~\lambda_{K}(\zz)-\lambda_{K}(\xx)=B^*\}; \qquad
\JJ:=\{k~:~ \lambda_{k}(\zz)-\lambda_{k}(\xx)=b^*\}.
\end{aligned}
\end{equation}
Then $b^*<B^*$ and $\I\cap \JJ=\emptyset$, but neither $\I$ nor $\JJ$  is empty, so that each set can have at most $n$ indices in them. Let us pick one $K\in \I$ and one $k\in \JJ$. Invoking Corollary \ref{l:pert_mkmi} we get  $\frac{\partial m_K}{\partial m_k}<0$, ($\frac{\partial \lambda_K}{\partial \lambda_k}<0$) hence we can obtain a new node system $\zz'$ in $X$, arbitrarily close to $\zz$, such that $\lambda_i(\zz)=\lambda_i(\zz')$ for all $i\ne k, K$, while the corresponding $\lambda_K(\zz')$ is a little smaller, than $\lambda_K(\zz)$ and simultaneously $\lambda_k(\zz')$ is a little greater than $\lambda_k(\zz)$. Note that for the indices $i\ne k,K$, the difference $\lambda_i(\zz')-\lambda_i(\xx)$ remained equal to $\lambda_i(\zz)-\lambda_i(\xx)$, as $\lambda_i$ was not changed.

Therefore, for the new node system $\zz'$ the respective values $b':=b(\zz') \ge b(\zz)$, and $B':=B(\zz') \le B(\zz)$, so that $\zz'\in W$, too, moreover, either $B'< B^*$, or at least $K$ is dropped from $\I$, and either $b'>b^*$, or at least $k$ is dropped from $\JJ$. Therefore, at most $n$ steps we  arrive at another node system $\zz'' \in X$, arbitrarily close to $\zz$, such that either $b(\zz'')>b(\zz)$ or $B(\zz'')<B(\zz)$. As the minimum is not decreased, and the maximum is not increased, the perturbed node system remains in the range where $\lambda_i(\zz'')\in [\lambda_i(\xx)+b^*,\lambda_i(\xx)+B^*]$. That is, $\zz, \zz', \zz'' \in Z$. However, $b^*$ was maximal and $B^*$ was minimal possible, so the very existence of $\zz''\in Z$ yields a contradiction. Consequently, the equality $b(\zz)=B(\zz)$ must hold, and the theorem is proved.
\end{proof}

Finally, let us present that the first part of the theorem entails in particular the positive answer to the Erd\H os Conjecture in our setting.

\begin{corollary}[Uniqueness and Sandwich Property]\label{cor:sandwich}
Let $w(t)=\exp(-t^2)$ and consider the respective Hermite interpolation. With the above notations the following holds.
If $\yy\in S$ is a minimax node system, then for any other node system
$\xx\in S$,
$\xx \ne \yy$ we necessarily have
$$
\mul(\xx):=
\min_{0\le i\le n+1}
m_i(\xx)< m(\yy) <\mol(\xx):=
\max_{0\le i\le n+1}
m_i(\xx).
$$
In particular, the extremal node system $\yy$ is unique.
\end{corollary}
\begin{proof} By Theorem \ref{th:minimax}, the the extremal node system $\yy$ is situated in $X$, and equioscillates. Therefore $\overline{m}(\yy)=m_0(\yy)=m_{n+1}(\yy)>1$. Now let $\xx \in X$ be a node system with $\mm(\xx)$ majorizing $\mm(\yy)$, or conversely, with $\mm(\yy)$ majorizing $\mm(\xx)$. Then according to Theorem \ref{th:nonmajorization} either $\xx=\yy$, which was excluded in the formulation of our Corollary, or we would have both $\xx, \yy \in X^c$, which contradicts to $\yy \in X$. Therefore, majorizing in any direction is impossible. The corollary is proved.
\end{proof}

\end{document}